\def\a{\alpha}

\def\b{\beta}

\def\CC{{\mathbb C}}
\def\cA{{\cal A}}

\def\cC{{\cal C}}

\def\cE{{\cal E}}

\def\cF{{\cal F}}
\def\cG{{\cal G}}

\def\cH{{\cal H}}

\def\cl{\colon}
\def\cL{{\cal L}}
\def\cM{{\cal M}}

\def\cO{{\cal O}}
\def\cod{\mbox{cod}}

\def\cQ{{\cal Q}}

\def\cU{{\cal U}}
\def\cV{{\cal V}}
\def\cW{{\cal W}}
\def\cX{{\cal X}}
\def\cY{{\cal Y}}

\def\d{\delta}

\def\deg{\text{\rom{deg}}}

\def\det{\text{\rom{det}}}

\def\e{\epsilon}

\def\Ext{\mbox{Ext}}
\def\es{\emptyset}

\def\g{\gamma}
\def\G{\Gamma}

\def\GG{{\mathbb G}}

\def\Hom{\mbox{Hom}}
\def\hra{\hookrightarrow}

\def\Id{\mbox{Id}}

\def\Im{\mbox{Im}}

\def\Ker{\mbox{Ker}}
\def\l{\lambda}
\def\la{\langle}

\def\L{\Lambda}
\def\lra{\longrightarrow}
\def\msk{\medskip}

\def\n{\noindent}

\def\NN{{\mathbb N}}

\def\om{\omega}

\def\op{\oplus}
\def\ot{\otimes}
\def\ov{\overline}
\def\O{\Omega}

\def\PP{{\mathbb P}}

\def\QQ{{\mathbb Q}}
\def\ra{\rangle}

\def\rk{\mbox{rk}}

\def\RR{{\mathbb R}}
\def\s{\sigma}

\def\Si{\Sigma}

\def\sm{\setminus}
\def\ss{\subset}

\def\T{\Theta}

\def\tm{\times}

\def\ul{\underline}

\def\uu{{\bf u}}

\def\vv{{\bf v}}

\def\wh{\widehat}
\def\wt{\widetilde}
\def\ww{{\bf w}}

\def\ZZ{{\mathbb Z}}

\documentclass[10pt,a4paper]{article}
\usepackage{amsthm}
\usepackage{amsmath}
\usepackage{amssymb}
\newtheorem{thm}{Theorem}[section]
\newtheorem{rmk}[thm]{Remark}
\newtheorem{dfn}[thm]{Definition}
\newtheorem{dfn-lmm}[thm]{Definition-Lemma}

\newtheorem{lmm}[thm]{Lemma}
\newtheorem{crl}[thm]{Corollary}
\newtheorem{prp}[thm]{Proposition}

\newtheorem{cnj}[thm]{Conjecture}
\newtheorem{hyp}[thm]{Hypothesis}
\newtheorem{clm}[thm]{Claim}

\newtheorem{stt}[thm]{Statement}

\begin{document}
 \title{Involutions and linear systems on holomorphic
 symplectic manifolds}
 \author{Kieran G. O'Grady\thanks{Supported by
 Cofinanziamento MURST 2002-03}\\
Universit\`a di Roma \lq\lq La Sapienza\rq\rq}
 \maketitle
 \section{Introduction}
A compact K\"ahler manifold is irreducible symplectic if it is
simply connected and it carries a holomorphic symplectic form
spanning the space of global holomorphic $2$-forms. A
$2$-dimensional irreducible symplectic manifolds is nothing else
but a $K3$ surface. The well established theory of periods of $K3$
surfaces has been a model for the theory in higher dimensions and
indeed Local Torelli~\cite{beau} and Surjectivity of the period
map~\cite{huy,huy-err} hold in any dimension. $K3$ surfaces are
remarkable not only for their periods: the complete linear system
associated to an ample divisor has very simple
behaviour~\cite{mayer} and furthermore one can describe
explicitely all $K3$ surfaces with an ample divisor whose
self-intersection is small. This paper deals with the question: do
similar properties hold for ample divisors on an irreducible
symplectic manifold of arbitrary dimension? Deformations of
$S^{[n]}$, the Hilbert scheme parametrizing length-$n$ subschemes
of a $K3$ surface $S$, are in many respects the symplest known
irreducible symplectic manifolds - recall that the generic such
deformation is not birational to $(K3)^{[n]}$ (Thm.~6, p.779
of~\cite{beau}) if $n\ge 2$. The L Conjecture~(\ref{lconj})
predicts that if $(X,H)$ is the generic couple with $X$ a
deformation of $(K3)^{[n]}$ and $H$ an ample divisor of square $2$
for Beauville's quadratic form then $|H|$ has no base-locus and
the map $X\to |H|^{\vee}$ has degree $2$ onto its image  and
furthermore $c_1(H)$ spans the subspace of $H^2(X)$ fixed by the
covering involution. When $n=1$ the conjecture is clearly true -
every degree-$2$ polarized $K3$ is a double cover of $\PP^2$. If
the L Conjecture is true then it follows that a deformation of
$(K3)^{[n]}$ carrying a divisor (not necessarily ample) of degree
$2$ has an anti-symplectic birational involution: this is a
non-trivial assertion, usually easier to test than the L
Conjecture - we call it the I Conjecture~(\ref{iconj}). Before
giving the precise statements we recall the properties of
Beauville's quadratic form (Thm.~5, p.~772 of~\cite{beau}). Let
$X$ be an irreducible symplectic manifold of (complex) dimension
$2n$: the quadratic form $(,)_X$ on $H^2(X)$ - a higher
dimensional analogue of the intersection form - is characterized
by the following properties:
 \begin{itemize}
\item[(1)] $(,)_X$ is integral indivisible non-degenerate, $(H^{p,q},H^{p',q'})_X=0$
if $p+p'\not=2$.
 \item[(2)] The signature of $(,)_X$ is $(3,b_2(X)-3)$. If $H$ is ample
 $(,)_X$ is positive
 definite on
 $$(H^{2,0}(X)\op H^{0,2}(X))_{\RR}\op\RR
 c_1(H).$$
 \item[(3)] There is a positive rational constant $c_X$
 such that the following formula of Fujiki holds (Thm.~(4.7)
 of~\cite{fuj}):
 \begin{equation}\label{fujformula}
 \int_X\a^{2n}=c_X \cdot (\a,\a)_X^n.
 \end{equation}
 \end{itemize}
Both $c_X$ and $(,)_X$ do not change if we modify the complex structure of $X$.
Beauville's form and the Fujiki constant of $(K3)^{[n]}$ are given in
Subsubsection~(\ref{hilbdescr}). Huybrechts (Lemma~(2.6) of~\cite{huy}) proved
that $(,)_X$ behaves well with respect to birational maps. More precisely let
$\phi\cl X\cdots> Y$ be a birational (i.e.~bimermomorphic) map between
irreducible symplectic manifolds and $H^2(\phi)\cl H^2(Y;\ZZ)\to H^2(X;\ZZ)$ be
defined by the K\"unneth decomposition of the Poincar\'e dual of the graph of
$\phi$.  Then $H^2(\phi)$ is an isometry of lattices. Furthermore, letting
$Bir(X)$ be the group of birational maps of $X$ to itself, the map
\begin{equation}\label{utile}
\begin{matrix} Bir(X) & \overset{H^2}{\lra} &
Isom(H^2(X;\ZZ),(,)_X)\\
\phi & \mapsto & H^2(\phi)
\end{matrix}
\end{equation}
is a homomorphism into the subgroup of integral Hodge isometries  of
$H^2(X;\CC)$. Whenever no confusion may arise we denote $(,)_X$ by $(,)$. For
$h\in H^{1,1}_{\ZZ}(X)$ with $(h,h)=2$ we let $R_h\cl H^2(X)\to H^2(X)$ be the
reflection in the span of $h$, i.e.
\begin{equation}\label{riflesso}
R_h(v)=-v+(v,h)h.
\end{equation}
\begin{dfn}\label{polsympl}
A couple $(X,H)$ is a {\sl degree-$k$ polarized irreducible symplectic variety}
if $X$ is an irreducible symplectic manifold and $H$ is an indivisible ample
divisor on $X$ with $(c_1(H),c_1(H))=k$.
\end{dfn}
Now let $(X,H)$ be a degree-$2$ polarized irreducible symplectic variety with
$X$ a deformation of $(K3)^{[n]}$. By a theorem of
Koll\'ar-Matsusaka~\cite{kolmat}  there is an $\ell_n>0$ depending only on $n$
such that the a priori rational map $X\cdots>|\ell_n H|^{\vee}$ is in fact a
regular embedding. Thus all  $(X,H)$ as above are realizable as subvarieties of
a $\PP^{d_n}$ with a fixed Hilbert polynomial $p_n$: let $\cQ_n$ be the Hilbert
scheme to which they \lq\lq belong\rq\rq and let $\cQ_n^0\ss\cQ_n$ be the open
subset given by
\begin{equation}\label{spiego}
\cQ_n^0:=\{t\in \cQ_n|\ \mbox{$X_t$ is irreducible symplecic}\},
\end{equation}
where $X_t$ is the subvariety of $\PP^{d_n}$ corresponding to $t\in\cQ_n$.
\begin{cnj}\label{lconj}{\rm[L Conjecture]}
Keep notation as above. There exists an open dense subset $U_n\ss\cQ_n^0$ such
that for $t\in U_n$ the following holds. Let $(X_t,H_t)$ be the degree-$2$
polarized irreducible symplectic variety corresponding to $t$. Then $|H_t|$ has
no base-locus and $f_t\cl X_t\to |H_t|^{\vee}$ is of degree $2$ onto its image
$Y_t$. In particular there exists an involution $\phi_t\cl X_t\to X_t$ such that
$f_t$ is the composition
$$X_t\overset{\pi_t}{\lra} X_t/\la\phi_t\ra\overset{\nu_t}{\lra}Y_t$$
where $\pi_t$ is the quotient map and $\nu_t$ is the normalization map. Let
$h_t:=c_1(H_t)$. Then
\begin{equation}\label{eccerifl}
H^2(\phi_t)=R_{h_t}.
\end{equation}
\end{cnj}
In the next section we will show  that if the above conjecture holds then also
 the following conjecture is true.
\begin{cnj}\label{iconj}{\rm[I Conjecture]}
Let $X$ be an irreducible symplectic manifold deformation equivalent to
$(K3)^{[n]}$. Suppose that $h\in H^{1,1}_{\ZZ}(X)$ and $(h,h)=2$. There exists a
birational involution $\phi\cl X\cdots>X$ such that for $\g\in H^2(X)$
\begin{equation}\label{azione}
H^2(\phi)(\g)=R_h(\g)-\sum_i (\g,\a_i)\b_i,
\end{equation}
where $\a_i\in H^{1,1}_{\QQ}(X)$ with $(\a_i,\cdot)$  equal to integration over
an effective analytic $1$-cycle and $\b_i\in H^{1,1}_{\ZZ}(X)$ is Poincar\'e
dual to an effective divisor.
\end{cnj}
A few comments. It follows from~(\ref{azione}) that $H^2(\phi)$
multiplies a symplectic form by $(-1)$: in particular $\phi$ is
not the identity! Since~(\ref{utile}) is a homomorphism
$H^2(\phi)$ is an involution: this imposes restrictions on the
$\a_i$'s and $\b_i$'s.

After proving that the L Conjecture implies the I Conjecture we will show - this
is easy - that the two conjectures are stable under deformations if certain
hypotheses are satisfied. More precisely: If $X,h,\phi$ are as in the I
Conjecture and furthermore $\phi$ is regular with $H^2(\phi)=R_h$ then $\phi$
extends to all small deformations of $X$ that keep $h$ of type $(1,1)$. If we
have $t_0\in\cQ_n^0$ such that $(X_{t_0},H_{t_0})$ behaves as stated in the L
Conjecture then the same holds for $(X_t,H_t)$ where $t$ varies in an open
subset of
 $\cQ_n^0$ containing (the orbit of)
$(X_{t_0},H_{t_0})$. In the next section we give examples of couples $X,h,\phi$
where $X$ is a deformation of $(K3)^{[n]}$, $h\in H^{1,1}_{\ZZ}(X)$ has degree
$2$ for Beauville's form and $\phi$ is a rational involution of $X$ with
$H^2(\phi)=R_h$. In our examples $X$ is always a moduli space of rank-$r$
torsion-free sheaves on a $K3$ surface $S$. The involutions were introduced by
Beauville~\cite{beau-aut} in the rank-$1$ case, by Mukai~\cite{mukaivb} when
$r\ge 2$. The new example in
 Subsubsection~(\ref{dualexample}) is a generalized Mukai
reflection. We spend some time proving that the action on $H^2$ is indeed the
reflection in  a class $h$ of square $2$: once this is proved we know that the
regular involutions are stable under small deformations of $(X,h)$. There are
examples of regular involutions in any (even) dimension as long as we allow $r$
to be arbitrarily large. We expect that all of the examples we give (with the
possible exception of the one in Subsubsection~(\ref{dualexample})) are \lq\lq
polarized\rq\rq deformation equivalent but we do not prove this, see
Section~(\ref{tuttiperuno}). One should notice that if we have $2$ regular
involutions $\phi_1,\phi_2$ on the same $X$ with $H^2(\phi_i)=R_{h_i}$ where
$h_1,h_2$ are independent then $\phi_1\circ\phi_2$ is an automorphism of
infinite order generating an interesting dynamical system. If furthermore $X$
and $\phi_1,\phi_2$ are defined over a number field $K$ one may study the action
of  $\phi_1\circ\phi_2$ on $X(\ov{K})$: this was done by Silverman~\cite{silver}
for $X$ a $K3$. We briefly discuss this in Subsection~(\ref{dueinv}).
Section~(\ref{exlinsyst}) is devoted to examples of degree-$2$ polarized $(X,H)$
where $X$ is a deformation of $(K3)^{[n]}$  and $|H|$ has the good behaviour
stated in the L Conjecture. We give examples in dimensions $4,6,8$. In doing so
we prove that the so-called Strange duality
statement~\cite{dontu,lepot,danila,lepot2} holds for certain couples of moduli
spaces of sheaves on a $K3$. We examine more closely the $4$-dimensional example
 (first given by Mukai~\cite{mukaisug}): a moduli space $X$ of rank-$2$ sheaves on a
$K3$ surface $S\ss\PP^6$ of degree-$10$, with $H$ a suitable ample divisor. We
have an identification $|H|^{\vee}\cong |I_S(2)|$ (Strange duality) and the
image of $X\to |I_S(2)|$ is the non-degenerate component, call it $Y$, of the
hypersurface parametrizing singular quadrics, the other component being a
hyperplane. $Y$ is a sextic $4$-fold in $\PP^5$, singular along a smooth
surface; Conjecture L predicts that the generic degree-$2$ polarized $(X,H)$
with $X$ a deformation of $(K3)^{[2]}$ is a double cover of a sextic $4$-fold in
$\PP^5$. This explains one of the main motivations for formulating our
conjectures. If the L Conjecture is true in dimension $4$ then we should have a
relatively explicit way of describing all degree-$2$ polarized $(X,H)$ with $X$
a deformation of $(K3)^{[2]}$ and hence also the relevant moduli space, call it
$M_2$. In this respect we notice that all known explicit constructions of
irreducible projective symplectic varieties give families of codimension $1$ in
the relevant moduli space, with one exception - the variety of lines on a cubic
$4$-fold~\cite{beaudon}: in this case we get a whole component of the moduli
space $M_6$ of degree-$6$ polarized $(X,H)$ with $X$ a deformation of
$(K3)^{[2]}$ and $(c_1(H),\cdot)$ a functional on $H^2(X;\ZZ)$ divisible by $2$.
If $M_6$ is irreducible then global Torelli for deformations of $(K3)^{[2]}$
follows from Voisin's Torelli Theorem~\cite{voitor} for cubic $4$-folds. We do
not know how to attack the problem of irreducibily of $M_6$; on the other hand
the L Conjecture should allow us to describe $M_2$, and once this is done we
should be in a better position to study the period map. Going back to Mukai's
example we notice that the dual hypersurface $Y^{\vee}\ss (\PP^5)^{\vee}$ is the
image of another symplectic variety $X^{\vee}$, in fact $X^{\vee}=S^{[2]}$, via
the complete linear system associated to a certain degree-$2$ divisor $H^{\vee}$
which is base-point free but not ample - it contracts a $\PP^2$. We expect that
there is an involution on the moduli space of degree-$2$ quasi-polarized $(X,H)$
($X$ a deformation of $(K3)^{[2]}$) which generalizes the above duality: we will
give some evidence for this in a forthcoming paper on the L Conjecture in
dimension $4$.
\section{The L Conjecture implies the I Conjecture}
\label{limpliesi}
 \setcounter{equation}{0}
Set $X_0=X$ and $h_0=h$. By Bogomolov's
Theorem~\cite{bog} $X_0$ has a smooth versal
deformation space and hence there exists a proper
submersive map $\pi\cl\cX\to B$ between manifolds
and a point $0\in B$ with the following
properties: $\pi^{-1}(0)\cong X_0$ and the
Kodaira-Spencer map
\begin{equation}\label{kodspe}
T_{B,0}\to H^1(T_{X_0})
\end{equation}
is an isomorphism. The germ $(\cX,X_0)\to (B,0)$ is identified
with the deformation space of $X_0$: we say that $\pi$ is a {\sl
representative of the versal deformation space of $X_0$}. We set
$X_t:=\pi^{-1}(t)$. The holomorphic symplectic form on $X_0$
defines an isomorphism $T_{X_0}\cong \O^1_{X_0}$ and hence
$H^1(T_{X_0})\cong H^1(\O^1_{X_0})$: since $h^1(\O^1_{X_0})=21$
(see Prop.~6, p.~768 of~\cite{beau}) we get by~(\ref{kodspe}) that
\begin{equation}\label{dimdef}
\dim B=21.
\end{equation}
We assume that $B$ is simply connected and
therefore there is a well-defined period map
$$\begin{matrix}
 B & \overset{P}{\lra} & \PP(H^2(X_0))\\
t & \mapsto & g_t^{-1}\left(H^{2,0}(X_t)\right)
 \end{matrix}$$
where $g_t\cl H^2(X_0)\to H^2(X_t)$ is given by parallel transport with respect
to the Gauss-Manin connection (see for example (9.2) of~\cite{voibook}) along
any path from $0$ to $t$. As is well-known (see Thm~(7.3), p.~254 of~\cite{bpv})
$P$ defines an isomorphism of a neighborhood of $0$ onto an open subset (in the
classical topology) of the smooth quadric
$$Q:=\{\ell\in\PP(H^2(X_0))|\ (\ell,\ell)=0\}.$$
Let $B(h_0)\ss B$ be the subset of $t$ such that
$g_t(h_0)\in H^2(X_t)$  is of type $(1,1)$: thus
$B(h_0)=P^{-1}(h_0^{\bot})$. Since $P$ is a local
isomorphism
\begin{equation}\label{smoothpol}
\mbox{$B(h_0)$ is smooth of codimension $1$ near $0$;}
\end{equation}
by~(\ref{dimdef}) we get that $\dim B(h_0)=20$. Let
$$B(h_0)_{am}:=\{t\in B(h_0)|\
\mbox{$h_t$ or $-h_t$ is the class of an ample
divisor}\}.$$
This is a Zariski-open subset of $B(h_0)$; we
claim that
\begin{equation}\label{ample}
B(h_0)_{am}\not=\es.
\end{equation}
In fact by the stated property of the period map
$P$ the set
$$B(h_0)_1:=\{t\in B(h_0)|\
H^{1,1}_{\ZZ}(X_t)=\ZZ h_t\}$$
is the complement of a countable union of proper
hypersurfaces of $B(h_0)$, hence there exists
$\ov{t}\in B(h_0)_1$; by Huybrechts' Projectivity
Criterion~\cite{huy,huy-err} $\ov{t}\in
B(h_0)_{am}$. Changing sign to $h_0$ if necessary
we can assume that $h_t$ is ample for $t\in
B(h_0)_{am}$. Now assume the L
Conjecture~(\ref{lconj}): then there is a Zariski
open non-empty subset
$$B(h_0)_{gd}\ss B(h_0)_{am}$$
such that for $t\in B(h_0)_{gd}$ the complete linear system $|H_t|$ enjoys the
properties stated in~(\ref{lconj}). (Here $H_t\in Pic(X_t)$ is the divisor class
such that $c_1(H_t)=h_t$.) Let $\phi_t\cl X_t\to X_t$ be the covering involution
and $\G_t\ss X_t\tm X_t$ be the graph of $\phi_t$. We prove the I
Conjecture~(\ref{iconj}) by considering the limit of $\G_t$ as $t\to 0$. We view
$\G_t$ as an element of the space $\cC_{2n}(\cX\tm_{B}\cX)$ parametrizing
effective compactly supported analytic $2n$-cycles constructed by
Barlet~\cite{barlet}. Let $\{t_k\}_{k\in\NN}$ be a sequence of points $t_k\in
B(h_0)_{gd}$ converging to $0$; such a sequence exists because $B(h_0)_{gd}$ is
a Zariski-open and dense subset of $B(h_0)$. Proceeding exactly as in the proof
of Theorem~(4.3) of~\cite{huy} we see that passing to a subsequence we can
assume that $\{\G_{t_k}\}_{k\in\NN}$ converges to an effective analytic
$2n$-cycle $\G_0$ on $X_0\tm X_0$. This (see the proof of Theorem~(4.3)
of~\cite{huy}) implies that there is a decomposition
$$\G_0=\G(\phi_0)+\sum_i n_i\O_i$$
where $\phi_0\cl X_0\cdots> X_0$ is a birational
map, $\G(\phi_0)$ is the graph of $\phi_0$,
$n_i>0$, $\O_i$ is irreducible and $\O_i\ss
D_i\tm E_i$ where $D_i,E_i\ss X_0$ are proper
subsets. Since $\G_{t_k}$ is invariant for the
involution of $X_{t_k}\tm X_{t_k}$ interchanging
the factors the same must hold for $\G_0$; this
implies that $\G(\phi_0)$ is invariant for the
involution, i.e.~$\phi_0$ is a birational
involution. Finally let's show that $H^2(\phi_0)$
is as stated in the I Conjecture~(\ref{iconj}).
For $\L$ an analytic cycle on $X_0\tm X_0$ we let
$H^2(\L)\cl H^2(X_0)\to H^2(X_0)$ be the map
defined by the K\"unneth component in
$H^{2}(X_0)\ot H^{2n-2}(X_0)$ of the Poincair\'e
dual of $\L$. Let $G_t\cl H^{*}(X_0\tm X_0)\to
H^{*}(X_t\tm X_t)$ be given by Gauss-Manin
parallel transport along any path going from $0$
to $t$: since $G_t([\G_0])=[\G_{t_k}]$ we have
\begin{equation}\label{gamact}
H^2(\G_0)=G_t^{-1}H^2(\G_{t_k})=G_t^{-1}H^2(\phi_{t_k})
=G_t^{-1}R_{h_t} =R_{g_t^{-1}h_t}=R_{h_0},
\end{equation}
where the third equality follows from~(\ref{eccerifl}). Now we determine
$H^2(\O_i)$. If $\cod(D_i,X)>1$ then $H^2(\O_i)=0$. Assume that $\cod(D_i,X)=1$
and let $C_i$ be a generic fiber of the map $\O_i\to D_i$ induced by the
projection $X_0\tm X_0\to X_0$ to the first factor; thus $C_i$ is a curve. Then
$$H^2(\O_i)(\g)=\left(\int_{p_{*}C_i}\g\right)[D_i]$$
where $p\cl X_0\tm X_0\to X_0$ is the projection to the second factor. This
equation together with~(\ref{gamact})  proves that~(\ref{azione}) holds.
\section{Stability results}\label{stabres}
 \setcounter{equation}{0}
Let $X_0$ be an irreducible symplectic manifold - not necessarily a deformation
of $(K3)^{[n]}$. A (regular) involution $\phi_0\cl X_0\to X_0$  is {\sl
anti-symplectic} if
\begin{equation}\label{meno}
 \phi_0^{*}\s_0=-\s_0
\end{equation}
where $\s_0$ is a holomorphic symplectic form on $X$. We have an orthogonal
decomposition into eigenspaces
$$H^2(X_0;\QQ)=
H^2(\phi_0)(+1)_{\QQ}\op_{\bot}H^2(\phi_0)(-1)_{\QQ}.$$
\subsection{Involutions}
 \setcounter{equation}{0}
Assume $X_0$ and $\phi_0$ are as above. By Bogomolov's Theorem~\cite{bog} $X_0$
has a smooth versal deformation space; let $\pi\cl\cX\to B$ be a representative
for the deformation space of $X_0$ and $X_t:=\pi^{-1}(t)$.  Let
\begin{equation}\label{extbir}
B(\phi_0):=\{t\in B|\ \mbox{$\exists \phi_t\cl
X_t\to X_t$ deformation of $\phi_0$}\}.
\end{equation}
Equation~(\ref{meno}) gives that $H^{2,0}(X_0)\ss H^2(\phi_0)(-1)$ and hence
$$L_{0}:=H^2(\phi_0)(+1)_{\QQ}\ss H^{1,1}_{\QQ}(X_0).$$
Let
$$B(L_0):=\{t\in B|\ g_t(L_0)\ss H^{1,1}_{\QQ}(X_t)\}$$
where $g_t\cl H^2(X_0)\to H^2(X_t)$ is given by
Gauss-Manin parallel transport along any path
connecting $0$ to $t$ (we assume that $B$ is
simply connected).
\begin{prp}\label{definv}
Keep notation and assumptions as above.  In a
neighborhood of $0$ we have $B(\phi_0)=B(L_0)$.
\end{prp}
\begin{proof}
By the universal property of the deformation space there exist an involution
$\tau\cl B\to B$ fixing $0$ and an isomorphism $\Phi\cl \tau^{*}\cX\to\cX$ of
families over $B$ whose restriction to $X_0$ is equal to $\phi_0$. (We allow
ourselves to shrink $B$.) Let $B^{\tau}\ss B$ be the fixed locus of $\tau$; then
$$B^{\tau}\ss B(\phi_0)\ss B(L_0).$$
Since $B^{\tau}$ is smooth it suffices to prove that
\begin{equation}\label{souguali}
T_{B^{\tau},0}=T_{B(L_0),0}.
\end{equation}
Let $\om_0$ be a symplectic form on $X_0$; contraction with $\om_0$ defines an
isomorphism
\begin{equation}\label{contratto}
T_{B,0}=H^1(T_{X_0})\overset{\sim}{\lra} H^1(\O^1_{X_0}).
\end{equation}
With this identification we have
$T_{B(L_0),0}=L_0^{\bot}=H^2(\phi_0)(-1)$. On the other hand
Isomorphism~(\ref{contratto}) gives an identification
$$T_{B^{\tau},0}=H^2(\phi_0)(-1)$$
because $\phi_0$ is anti-symplectic. This proves~(\ref{souguali}).
\end{proof}
As an immediate consequence we have the following
result.
\begin{crl}\label{stabinv}
Keep notation and hypotheses as above. Assume furthermore that
$H^2(\phi_0)=R_{h_0}$ where $h_0\in H^{1,1}_{\ZZ}(X_0)$. Then $\phi_0$ extends
to all small deformations of $X_0$ that keep $h_0$ of type $(1,1)$.
\end{crl}
\subsection{Linear systems}
 \setcounter{equation}{0}
Let $(X_0,H_0)$ be a degree-$2$ polarized deformation of $(K3)^{[n]}$. Let
$\cQ_n^0$ be the open subset of a Hilbert scheme given by~(\ref{spiego}); thus
we may think that $0\in \cQ_n^0$ and that $X_0$ is the (embedded) variety
corresponding to $0$, with $\cO_{X_0}(1)\cong\cO_{X_0}(\ell_n H_0)$. Notice that
$\cQ_n^0$ is smooth at $0$ because by~(\ref{smoothpol}) the deformation space of
$(X_0,H_0)$ is smooth: let $\cA\ss\cQ_n^0$ be the irreducible component
containing $0$.
\begin{prp}\label{stablin}
Keep notation as above. Suppose that $|H_0|$ has no base-locus and that $f_0\cl
X_0\to |H_0|^{\vee}$ is of degree $2$ onto its image $Y_0$, hence there exists
an involution $\phi_0\cl X_0\to X_0$ such that $f_0$ is the composition
$$X_0\overset{\pi_0}{\lra}X_0/\la\phi_0\ra\overset{\nu_0}{\lra}Y_0,$$
where $\pi_0$ is the quotient map and $\nu_0$ is the normalization map. Suppose
also that
\begin{equation}\label{reflex}
H^2(\phi_0)=R_{h_0}
\end{equation}
where $h_0:=c_1(H_0)$. Then there exists an open non-empty subset $\cV\ss\cA$
such that the statements above hold when we replace $(X_0,H_0)$ by $(X_t,H_t)$.
\end{prp}
\begin{proof}  Having no base-locus is an open property; since $|H_0|$ has no
base-locus we get that $|H_t|$ has have no base-locus for $t$ varying in an open
non-empty subset of $\cA$. Now consider the other statements. First locally
around $0$ we can extract the $\ell_n$-th root of $\cO_{X_t}(1)$. More precisely
 there exist a quasi-projective manifold $\cU$, a finite map $\rho\cl \cU\to\cA$
and a point $\wt{0}\in \cU$ with the following properties: $\rho(\wt{0})=0$,
$\rho$ is submersive at $\wt{0}$ and the pull-back $\zeta\cl\cX\to \cU$ of the
tautological family over $\cA$ carries a divisor class $\cH$ such that for $t\in
U$ we have $(c_1(H_t),c_1(H_t))=2$ and $\cO_{X_t}(1)\cong \cO_{X_t}(\ell_n
H_t)$. (Here $X_t:=\pi^{-1}(t)$ and $H_t:=\cH|_{X_t}$.) By
Corollary~(\ref{stabinv}) there exists a Zariski-open $\cU_{inv}\ss \cU$ such
that for $t\in \cU_{inv}$ we have an involution $\phi_t\cl X_t\to X_t$ with
$H^2(\phi_t)=R_{h_t}$. Let $\cX_{inv}:=\pi^{-1}(\cU_{inv})$; then we have an
involution $\Phi\cl\cX_{inv}\to\cX_{inv}$ restricting to $\phi_t$ on each $X_t$.
Let $\wt{\cY}:=\cX_{inv}/\langle\Phi\rangle$ be the quotient. The map
$\xi\cl\wt{\cY}\to \cU_{inv}$ is analytically locally trivial hence $\xi$ is a
flat family. Let $\cH_{inv}$ be the restriction of $\cH$ to $\cX_{inv}$. We
claim that $\cH_{inv}$ descends to a divisor class $\ov{\cH}$ on $\wt{\cY}$:
this follows at once from the fact that by hypothesis $H_0$ descends to the
divisor class on $\wt{Y}_0$ given by $\nu_0^{*}\cO_{Y_0}(1)$. Letting $f_t\cl
X_t\to \wt{Y}_t$ be the quotient map we have the pull-back
\begin{equation}\label{tira}
f_t^{*}\cl H^0(\wt{Y}_t;\ov{H}_t)\to H^0(X_t;H_t).
\end{equation}
We claim that $f_t^{*}$ is an isomorphism for all $t\in \cU_{inv}$. Since
$f_t^{*}$ is injective it suffices to check that
\begin{equation}\label{equalh0}
h^0(\wt{Y}_t;\ov{H}_t)= h^0(X_t;H_t).
 \end{equation}
We claim that
\begin{equation}\label{vanishing}
h^p(\wt{Y}_t;\ov{H}_t)=0=h^p(X_t;H_t),\quad p>0.
\end{equation}
This follows from the classical Kodaira vanishing for $X_t$ because $K_{X_t}\sim
0$ and for $\wt{Y}_t$ we apply for example Theorem~(1-2-5) of~\cite{kmm}; notice
that $\wt{Y}_t$ has terminal singularities and $K_{\wt{Y}_t}\equiv 0$ (for this
we need $\dim X_0\ge 4$, if $\dim X_0=2$ we are considering $K3$ surfaces and
the proposition is trivially verified). From~(\ref{vanishing}) we get that in
order to prove~(\ref{equalh0}) it suffices to show that
$\chi(\wt{Y}_t;\ov{H}_t)= \chi(X_t;H_t)$. Since $\cX_{inv}\to \cU_{inv}$ and
$\wt{\cY}\to \cU_{inv}$ are flat families and $\cU_{inv}$ is connected it is
enough to check that $\chi(\wt{Y}_0;\ov{H}_0)= \chi(X_0;H_0)$.
By~(\ref{vanishing}) this is equivalent to $h^0(\wt{Y}_0;\ov{H}_0)=
h^0(X_0;H_0)$: this we know by hypothesis. We have proved that~(\ref{tira}) is
an isomorphism. It follows that there is an open non-empty $\cU'_{inv}\ss
\cU_{inv}$ such that for $t\in \cU'_{inv}$ the  map $X_t\to Y_t\ss |H_t|^{\vee}$
factors through the quotient map $X_t\to \wt{Y}_t$ and that the induced map
$\wt{Y}_t\to Y_t$ is the normalization map. Finally $H^2(\phi_t)$ is constant
with respect to the Gauss-manin connection because $\phi_t$ is regular for all
$t\in \cU'_{inv}$ and hence $H^2(\phi_t)=R_{h_t}$ because~(\ref{reflex}) holds.
 \end{proof}

\section{Examples: involutions}\label{exinv}
We give examples of birational involutions on deformations of $(K3)^{[n]}$ whose
action on $H^2$ is the reflection in a $(1,1)$-class $h$. In many cases the
involution is regular and hence by Corollary~(\ref{stabinv}) it will extend to
all small deformations keeping $h$ of type $(1,1)$.
\subsection{Beauville's examples}\label{beauvex}
 \setcounter{equation}{0}
These are involutions of $S^{[n]}$ where $S$ is a $K3$ surface; they were
introduced by Beauville, see pp.~20-25 of~\cite{beau-aut}.
\subsubsection{Hilbert scheme of points on a $K3$}\label{hilbdescr}
We recall the description of $H^2(S^{[n]})$ and
its Beauville form (see Prop.~6, p.768
of~\cite{beau}, the remark following it and
pp.~777-778). Let $S^{(n)}$ be the symmetric
product of $n$ copies of $S$; we have a natural
\lq\lq symmetrization\rq\rq map $s\cl
H^2(S;\ZZ)\to H^2(S^{(n)};\ZZ)$. The cycle map
$c\cl S^{[n]}\to S^{(n)}$ gives rise to $c^{*}\cl
H^2(S^{(n)};\ZZ)\to H^2(S^{[n]})$. Composing $s$
and $c^{*}$ we get the map
\begin{equation}\label{mumap}
\mu\cl H^2(S;\ZZ)\to H^2(S^{[n]};\ZZ)
\end{equation}
which is an injection onto a saturated subgroup
of $H^2(S^{[n]};\ZZ)$. Let
\begin{equation}\label{nonred}
\Delta_n:=\{[Z]\in S^{[n]}|\ \mbox{$Z$ is
non-reduced}\},
\end{equation}
i.e.~the exceptional divisor of $c$: there exists
a (unique) divisor-class $\Xi_n$ such that
$2\Xi_n\sim\Delta_n$, set $\xi_n:=c_1(\Xi_n)$.
There is a direct sum decomposition
\begin{equation}\label{h2hilb}
H^2(S^{[n]};\ZZ)=\mu
\left(H^2(S;\ZZ)\right)\op_{\bot}\ZZ\xi_n
\end{equation}
orthogonal with respect to Beauville's form, and
furthermore
$$\begin{array}{rcl}
(\mu(\a),\mu(\b)) & =&\int_S\a\wedge\beta\\
(\xi_n,\xi_n) & =&-2(n-1).
\end{array}$$
Since $\mu$ is a morphism of Hodge structures and $\xi_n$ is of type $(1,1)$
Equality~(\ref{h2hilb}) determines the Hodge structure of $H^2(S^{[n]})$.
Finally we recall that the Fujiki constant (see~(\ref{fujformula})) of
$X=S^{[n]}$ is given by
\begin{equation}\label{fujconst}
c_X=\frac{(2n)!}{n! 2^n}\,.
\end{equation}
\subsubsection{The involution}
Assume that $D_{2g-2}$ is a globally generated
ample divisor on $S$ with
$$D_{2g-2}\cdot D_{2g-2}=2g-2.$$
 If $[Z]\in S^{[g-1]}$ is generic
then $|I_Z(D_{2g-2})|$ is one-dimensional and its
base-locus (as a linear system of divisors on
$S$) equals $Z\coprod W$, where $W$ is a
length-$(g-1)$ subscheme of $S$: Beauville
defines a birational involution
\begin{equation}\label{beauvinv}
\phi\cl S^{[g-1]}\cdots> S^{[g-1]}
\end{equation}
by setting $\phi([Z]):=[W]$ for the generic
$[Z]\in S^{[g-1]}$. For $g=2$ the map $\phi$ is
the involution defined on a $K3$ of degree-$2$.
For $g\ge 3$ and $D_{2g-2}$ very ample we have
$S\ss \PP^{g}$ and for $[Z]$ generic $Z\coprod
W=\langle Z\rangle\cap S$ where $\langle
Z\rangle$ is the $(g-2)$-dimensional span of $Z$.
We will study $H^2(\phi)$. Let
\begin{equation}\label{hgdef}
h_{g}:=
 \left(\mu\left(c_1(D_{2g-2})\right)-\xi_{g-1}\right)\in
 H^{1,1}_{\ZZ}(S^{[g-1]}).
 \end{equation}
Notice that $(h_g,h_g)=2$.
\begin{prp}\label{beauvact}
For $(S,D_{2g-2})$ varying in an open dense
subset of the relevant moduli space of polarized
$K3$ surfaces $H^2(\phi)$ equals the reflection
in the span of $h_g$ i.e.
\begin{equation}\label{reflhg}
H^2(\phi)(\g)=R_{h_g}:= -\g+(\g,h_g)h_g.
\end{equation}
\end{prp}
\begin{rmk}
The \lq\lq relevant moduli space...\rq\rq means
the following. Let $D_{2g-2}\sim k D_{2\ov{g}-2}$
where $k\in\NN$ and $c_1(D_{2\ov{g}-2})$ is
indivisible: the relevant moduli space is that of
degree-$(2\ov{g}-2)$ polarized $K3$'s.
\end{rmk}
\n
 {\bf Proof of Proposition~(\ref{beauvact}).}
First consider the case $g=2$: then
$S/\langle\phi\rangle\cong\PP^2$ hence the
$(+1)$-eigenspace of $H^2(\phi)$ is generated by
$h_2$ and thus $H^2(\phi)$ is the reflection in
the span of $h_2$.  Now assume that $g\ge 3$. By
general considerations there is an open dense
subset $U$ of the moduli space such that
$H^2(\phi)$ is \lq\lq constant\rq\rq over $U$.
Shrinking $U$ if necessary we can assume that
$D_{2g-2}$ is very ample for all
$[(S,D_{2g-2})]\in U$. For $[(S,D_{2g-2})]\in U$
let
$$f\cl S^{[g-1]} \cdots>  {\bf Gr}(1,|D_{2g-2}|)$$
 be the rational map sending
$[Z]$ to $|I_Z(D_{2g-2})|$, and let $L$ be the
(very ample) line bundle on ${\bf
Gr}(1,|D_{2g-2}|)$ defined by the Pl\"ucker
embedding. Since $\phi$ commutes with $f$ we have
\begin{equation}
\phi^{*}f^{*}c_1(L)=f^{*}c_1(L).\label{inv}
\end{equation}
We claim that
\begin{equation}
f^{*}c_1(L)=h_g.\label{pull}
\end{equation}
It suffices to prove~(\ref{pull}) for one
$(S,D_{2g-2})\in U$ because $U$ is irreducible.
By Hodge theory there exists $(S,D_{2g-2})\in U$
such that
\begin{equation}\label{noprim}
H^{1,1}_{\bf \QQ}(S)=\QQ c_1(D_{2g-2}).
\end{equation}
By~(\ref{h2hilb}) we have
\begin{equation}\label{nersev}
 H^{1,1}_{\bf
\QQ}(S^{[g-1]})=\QQ\mu\left(c_1(D_{2g-2})\right)
\op\QQ\xi_{g-1}.
 \end{equation}
Thus $f^{*}c_1(L)=x\mu\left(c_1(D_{2g-2})\right)+y\xi_{g-1}$. We get $x=1$ and
$y=-1$ by intersecting with the algebraic $1$-cycles
$$\begin{array}{rcl}
\G & := & \{[p_1+\cdots p_{g-2}+p]|\ p\in C\},\\
\L & := & \{[p_1+\cdots p_{g-3}+Z']|\ Z' \mbox{
non-reduced}\},
\end{array}$$
where $p_1,\ldots,p_{g-2}$ are fixed and $C\in
|D_{2g-2}|$; one must recall that
\begin{equation}
\langle
c_1(\Delta_{g-1}),\L\rangle=-2.\label{negint}
\end{equation}
From~(\ref{inv})-(\ref{pull}) we get that
\begin{equation}
\phi^{*}h_g=h_g.\label{hfixed}
\end{equation}
Now we determine the action of $\phi^{*}$ on the remaining part of
$H^2(S^{[g-1]})$. It suffices to prove that~(\ref{reflhg}) holds for
$(S,D_{2g-2})\in U$ such that~(\ref{noprim}) holds. Since~(\ref{utile}) is a
homomorphiam $H^2(\phi)$ is an isometric involution; by~(\ref{nersev})
and~(\ref{hfixed}) we get that the restriction of $H^2(\phi)$ to $H^{1,1}_{\bf
\QQ}(S^{[g-1]})$ is either  the identity or the reflection $\QQ{h_g}$. To show
that the latter holds it suffices to check that
\begin{equation}\label{noteq}
\phi^{*}c_1(\Delta_{g-1})\not= c_1(\Delta_{g-1}).
 \end{equation}
By~(\ref{negint}) any effective divisor
homologous to $\Delta_{g-1}$ must be equal to
$\Delta_{g-1}$; since
$\phi^{*}(\Delta_{g-1})\not=\Delta_{g-1}$ we
get~(\ref{noteq}). Now consider
$T(S^{[g-1]}):=H^{1,1}_{\bf
\QQ}(S^{[g-1]})^{\bot}$: it is left invariant by
the isometry $H^2(\phi)$.  To finish the proof it
will suffice to show that the restriction of
$H^2(\phi)$ to $T(S^{[g-1]})$ is equal to $(-1)$.
Since the eigenspaces of the restriction of
$H^2(\phi)$ to $T(S^{[g-1]})$ are Hodge
substructures and  $T(S^{[g-1]})$ has no
non-trivial Hodge substructures, the restriction
of $H^2(\phi)$ to $T(S^{[g-1]})$ is equal to $\pm
1$. Thus it suffices to show that
\begin{equation}
\phi^{*}(\s^{[g-1]})=-\s^{[g-1]}\label{minus}
\end{equation}
where $\s^{[g-1]}$ is the symplectic form on
$S^{[g-1]}$ induced by a symplectic form $\s$ on
$S$ (see~\cite{beau}, p.~766). Letting
$\phi([Z])=[W]$ we have
$$Z+W\sim D_{2g-2}\cdot D_{2g-2},$$
and hence Equality~(\ref{minus}) follows from Mumford's Theorem on $0$-cycles
(see Prop~(22.24) of~\cite{voibook}).
 \qed
\subsubsection{More on the involution}\label{specifico}
 Let
 $f_S\cl S\to |D_{2g-2}|^{\vee}$ be the natural
 map and let $H_g$ be the divisor class
 on $S^{[g-1]}$ such that $c_1(H_g)=h_g$.
 Suppose first that $g=3$, and consider the three possible
  cases
 \begin{itemize}
 \item[(a)]
 $D_4$ is very ample and $Im(f_S)$
does not contain lines,
 \item[(b)]
$D_4$ is very ample and $Im(f_S)$ does contain
lines $\ell_1,\ldots,\ell_k$,
 \item[(c)]
 $D_4$ is not very ample, i.e.~$f_S$ is
 $2$-to-$1$ onto a quadric.
 \end{itemize}
In case~(a) the divisor class $H_3$ is  ample. Furthermore $\phi$ is regular
(Proposition~(11) of~\cite{beau-aut}). Thus $H^2(\phi)$ is constant on the open
subset parametrizing $(S,D_{2g-2})$ for which~(a) holds; by
Proposition~(\ref{beauvact}) we get that $H^2(\phi)$ is the reflection in $\ZZ
h_3$. Applying Corollary~(\ref{stabinv}) we get that $\phi$ extends to all small
deformations of $S^{[2]}$ keeping $h_3$ of type $(1,1)$ - notice that the
generic such deformation is not of the type $(K3)^{[2]}$. In case~(b) the
divisor class $H_3$ is globally generated and big but not ample. Furthermore
$\phi$ is not regular (Proposition~(11) of~\cite{beau-aut}). Beauville shows
that to resolve the indeterminacies of $\phi$ it suffices to blow up
$\ell_1^{(2)}\cup\cdots\cup\ell_k^{(2)}$, and in fact $\phi$ lifts to a regular
involution on the blow-up: thus $\phi$ is the flop of
$\ell_1^{(2)}\cup\cdots\cup\ell_k^{(2)}$. As is easily checked $H^2(\phi)$ is
the reflection in $\ZZ h_3$. In case~(c) the map $\phi$ is not regular, in
particular $H_3$ is not ample. Furthermore $H^2(\phi)$ is not the reflection in
$\ZZ h_3$. Now consider $g\ge 4$: then the map $\phi$ is never regular (p.~24
of~\cite{beau-aut}), in particular $H_g$ is not ample. If $g=4,5$ and $S$ is
generic then Beauville (p.~25 of~\cite{beau-aut}) shows that the indeterminacies
of $\phi$ are resolved by a single blow-up with smooth center ($\phi$ is a Mukai
elementary modification~\cite{mukaisympl}).
\subsection{Mukai reflections}
 \setcounter{equation}{0}
These are involutions of moduli spaces of sheaves on a $K3$ surface. Beauville's
involutions are never regular in dimension greater than $4$: on the other hand
Mukai reflections give examples of regular involutions on deformations of
$(K3)^{[n]}$ for arbitrary $n$. The action of a Mukai reflection on $H^2$ is the
reflection in a class of square $2$.
\subsubsection{Moduli of sheaves on a $K3$ surface $S$}
We recall basic definitions and results. Let $F$ be a sheaf on $S$; following
Mukai~\cite{mukaivb} one sets
$$v(F):=ch(F)\sqrt{Td(S)}=ch(F)(1+\eta)\in H^{*}(S;\ZZ),$$
where $\eta\in H^4(S;\ZZ)$ is the orientation class. For
 $\a\in H^{*}(S)$ with degree-$q$ component given by
 $\a_{q}$ we
set
\begin{equation}
\a^{\vee}:=\a_0-\a_2+\a_4.\label{star}
\end{equation}
On $H^{*}(S)$ we have {\sl Mukai's bilinear symmetric form\/} defined by
\begin{equation}\label{mukform}
\langle u,w\rangle:=-\int_S u\wedge w^{\vee}.
\end{equation}
By Hirzebruch-Riemann-Roch we have
\begin{equation}\label{eulchar}
 \langle
v(E),v(F)\rangle=-\chi(E,F):=
-\sum\limits_{i=0}^{2}(-1)^i\dim
\mbox{Ext}^i(E,F).
\end{equation}
Let
\begin{equation}\label{explicitv}
 \vv=r+\ell+s\eta\in H^0(S;\ZZ)_{\ge 1}\op H^{1,1}_{\ZZ}(S)\op H^4(S;\ZZ).
 \end{equation}
Given an ample divisor $D$ on $S$ we let $\cM(\vv)$ be the moduli space of
Gieseker-Maruyama $D$-semistable torsion-free sheaves $F$ on $S$ with
$v(F)=\vv$; this is a projective variety~\cite{gieseker,maruyama}. An example:
let $D_{2g-2}$ be a divisor on $S$ with $D_{2g-2}\cdot D_{2g-2}=2g-2$ and let
\begin{equation}\label{hilbvect}
 \vv:=1+c_1(D_{2g-2})+\eta.
  \end{equation}
Then $\cM(\vv)$ parametrizes sheaves $I_Z(D_{2g-2})$ where $Z$ is a
length-$(g-1)$ subscheme of $S$, hence $\cM(\vv)=S^{[g-1]}$. For the sake of
simplicity we omit reference to $D$ in the notation for $\cM(\vv)$; however one
should keep in mind that if we change $D$ the moduli space $\cM(\vv)$ might
change.  Assume that $F$ is stable and that $v(F)=\vv$: the tangent space of
$\cM(\vv)$ at the point $[F]$ corresponding to $F$ is canonically identified
with $\mbox{Ext}^1(F,F)$ (see Cor~(4.5.2) of~\cite{huy-lehn}). Stability of $F$
and Serre duality give
$$1=\dim\mbox{Hom}(F,F)=\dim\mbox{Ext}^2(F,F),$$
hence~(\ref{eulchar}) gives
\begin{equation}\label{tandim}
\dim\mbox{Ext}^1(F,F)=2+\langle \vv,\vv\rangle.
\end{equation}
By a theorem of Mukai~\cite{mukaivb} we know that $\cM(\vv)$ is smooth at $[F]$
and thus
\begin{equation}\label{dimmod}
 \dim_{[F]}\cM(\vv)=2+\langle
\vv,\vv\rangle\quad \mbox{if $F$ is stable.}
 \end{equation}
It has been proved that under certain hypotheses on $D$ and $\vv$ the moduli
space $\cM(\vv)$ is an irreducible symplectic variety deformation equivalent to
$S^{[n]}$ where $2n=2+\langle \vv,\vv\rangle$. In order to state a result which
suffices for our purposes we give a definition.
\begin{dfn}\label{vgeneric}
Keep notation and assumptions as above. The ample divisor $D$ is $\vv$-generic
if there exists no couple $(r_0,\ell_0)$ consisting of an integer $0<r_0<r$ and
 $\ell_0\in H^{1,1}_{\ZZ}(S)$ such that
$$\left(r\ell_0-r_0 \ell\right)\cdot D=0,\quad
-(r^2\langle\vv,\vv\rangle+2r^4)\le
 4\left(r\ell_0-r_0 \ell\right)^2<0.$$
\end{dfn}
Given $\vv$ there exists a $\vv$-generic ample $D$~\cite{ogradyvb}. The
following lemma is proved in~\cite{ogradyvb}.
\begin{lmm}\label{vgenlem}
Keep notation and assumptions as above. Assume that $r+\ell$ is indivisible and
that $D$ is $\vv$-generic. Let  $F$ be a $D$-slope semistable sheaf with
$v(F)=\vv$; then $F$ is slope-stable. In particular $\cM(\vv)$ is smooth.
\end{lmm}
The following theorem was proved by Yoshioka~\cite{yoshi} (see~\cite{ogradyvb}
for the case when $\ell$ is indivisible).
\begin{thm}\label{yoshi}{\rm[Yoshioka]}
Keep notation and assumptions as above. Assume that $r+\ell$ is indivisible and
that $D$ is $\vv$-generic. Then $\cM(\vv)$ is an irreducible symplectic variety
deformation equivalent to $S^{[n]}$ where $2n=2+\langle \vv,\vv\rangle$.
\end{thm}
Under these hypotheses there is a beautiful description of $H^2(\cM(\vv))$ and
its Beauville form given by Mukai~\cite{mukaisug} and proved by
Yoshioka~\cite{yoshi} (see~\cite{ogradyvb} for the case when $\ell$ is
indivisible). First we need some preliminaries (see~\cite{mukaivb}). A {\sl
quasi-family of sheaves on $S$ parametrized by $T$ with Mukai vector $\vv$}
consists of a sheaf $\cF$ on $S\tm T$ flat over $T$ with $\cF|_{S\tm\{t\}}\cong
F^{\op d}$ where $v(F)=\vv$ and $d$ is some positive integer independent of $t$;
we set $\s(\cF):=d$. Two quasi-families $\cF,\cG$ on $S$ parametrized by $T$
with Mukai vector $\vv$ are {\sl equivalent} if there exist vector-bundles
$\cV,\cW$ on $T$ such that $\cF\ot p_T^{*}\cV\cong \cG\ot p_T^{*}\cW$, where
$p_T\cl S\tm T\to T$ is the projection. Given a quasi-family $\cF$ as above we
let $\theta_{\cF}\cl H^{*}(S)\to H^2(T)$ be defined by
\begin{equation}\label{thetamap}
\theta_{\cF}(\a):= \frac{1}{\s(\cF)} p_{T,*} \left[ch(\cF)(1+p_S^{*}\eta)
p_S^{*}(\a^{\vee})\right]_6
\end{equation}
where $p_S,p_T\cl S\tm T\to S,T$ are the projections  and $[\cdots]_q$ denotes
the component of $[\cdots]$ in $H^q(S\tm T)$. An easy computation gives the
following result.
\begin{lmm}\label{nondipende}
Keeping notation as above assume that $\cF,\cG$
are two equivalent quasi-families of sheaves on
$S$ parametrized by $T$ with Mukai vector $\vv$.
If $\a\in\vv^{\bot}$ then
$$\theta_{\cF}(\a)=\theta_{\cG}(\a).$$
\end{lmm}
Mukai (Thm.~(A.5) of~\cite{mukaivb}) showed that there exists a {\sl
tautological quasi-family $\cE$ on $S$ parametrized by $\cM(\vv)$\/} i.e.~a
quasi-family $\cE$ of sheaves on $S$ parametrized by $\cM(\vv)$ with Mukai
vector $\vv$ such that $\cE|_{S\tm [F]}\cong F^{\op \s(\cF)}$. Furthermore the
proof of Thm.~(A.5) of~\cite{mukaivb} shows that any two tautological
quasi-familes are equivalent. Thus by Lemma~(\ref{nondipende}) we get a
well-defined linear map
\begin{equation}\label{thetaiso}
 \theta_{\vv}\cl\vv^{\bot}\lra H^2(\cM(\vv))
 \end{equation}
by setting $\theta_{\vv}:=\frac{1}{\s(\cF)}\theta_{\cF}$ where $\cF$ is any
tautological quasi-family of sheaves on $S$ parametrized by $\cM(\vv)$. Now
define a weight-two Hodge structure on $H^{*}(S)$ by setting $F^0:=H^{*}(S)$,
$F^1:=H^0(S)\op F^1 H^2(S)\op H^4(S)$ and $F^2:=F^2 H^2(S)$. Since $\vv$ is
integral of type $(1,1)$ the orthogonal $\vv^{\bot}$ inherits a lattice
structure and a Hodge structure from $H^{*}(S)$. The following result was proved
by Yoshioka~\cite{yoshi} (see~\cite{ogradyvb} for the case when $\ell$ is
primitive and~\cite{mukaivb} for $\vv$ isotropic).
\begin{thm}\label{h2iso}{\rm[Yoshioka]}
Keep notation and assumptionms as above. Suppose that $r+\ell$ is indivisible,
 that $\langle \vv,\vv\rangle\ge 2$ and that $D$ is $\vv$-generic.
 Then $\theta_{\vv}$ is an isomorphism of integral
Hodge structures and defines an isometry between $\vv^{\bot}$ and
$\left(H^2(\cM(\vv);\ZZ),(\cdot,\cdot)\right)$.
\end{thm}
An example: if $\vv$ is given by~(\ref{hilbvect}) we get~(\ref{h2hilb}) with
$n=g-1$.
\subsubsection{Definition of Mukai reflections}\label{mukairefl}
Set $r=s$ in~(\ref{explicitv}) i.e.
\begin{equation}\label{vsymmetric}
 \vv=r+\ell+r\eta, \quad r\ge 1.
  \end{equation}
Under certain hypotheses there exists a {\sl Mukai reflection\/} on $\cM(\vv)$,
i.e.~a birational involution
\begin{equation}\label{phiv}
 \phi_{\vv}\cl
\cM(\vv)\cdots>\cM(\vv).
 \end{equation}
We will make the following assumption.
\begin{hyp}\label{picint}
If $A$ is a divisor on $S$ the intersection
number $A\cdot D$ is a multiple of $\ell\cdot D$.
 \end{hyp}
\begin{rmk}\label{itisgeneric}
If Hypothesis~(\ref{picint}) holds then $\ell$ is
indivisible and $D$ is $\vv$-generic. By
Theorem~(\ref{yoshi}) we get that $\cM(\vv)$ is a
deformation of $S^{[n]}$ where $2n=\langle
\vv,\vv,\rangle+2$. If furthermore
$\langle\vv,\vv\rangle\ge 2$ then by
Theorem~(\ref{picint})  the Hodge and lattice
structures on $H^2(\cM(\vv))$ are isomorphic to
those of $\vv^{\bot}$.
\end{rmk}
We also add the following assumption:
\begin{equation}\label{eldpos}
\ell\cdot D>0.
\end{equation}
Let $[F]\in\cM(\vv)$: then
\begin{equation}\label{h2andchi}
h^2(F)=0,\quad \chi(F)=\chi(\cO_S,F)=2r.
\end{equation}
To get the first equality notice that $H^2(F)\cong \Hom(F,\cO_S)^{\vee}$
 by Serre duality and that by slope-semitability of $F$ we have
 $\Hom(F,\cO_S)=0$. The second equality follows from~(\ref{eulchar}). From~(\ref{h2andchi}) we get
that
\begin{equation}\label{lowbound}
\mbox{$h^0(F)\ge 2r$ for $[F]\in\cM(\vv)$.}
\end{equation}
Let
 \begin{equation}\label{efgamma}
\wt{F}:=\Im(H^0(F)\ot\cO_S\to F)
 \end{equation}
be the subsehaf of $F$ generated by global
sections.
\begin{lmm}\label{stablem}{\rm[Markman]}
Keep notation and hypotheses as above. Then
\begin{itemize}
\item[(1)]
$F/\wt{F}$ is Artinian.
\item[(2)]
 the sheaf $E$ fitting into the exact sequence
\begin{equation}\label{edef}
0\to E\to H^0(F)\ot\cO_S\to \wt{F}\to 0
\end{equation}
is locally-free and slope-stable.
\end{itemize}
\end{lmm}
\begin{proof} (1) follows
 from Lemma~(3.5), p.~661 in~\cite{markman}. (2): $E$ is locally-free
 because $\wt{F}$ is a torsion-free sheaf on the smooth
 surface $S$ and hence its projective dimension
 is at most $1$. That $E$ is slope-stable is
 proved in~\cite{markman}, pp.~682-684 \lq\lq The
 case $a+b+t>a$\rq\rq.
 \end{proof}

Now let $U(\vv)\ss\cM(\vv)$ be the subset defined
by
\begin{equation}\label{uvdef}
U(\vv):=\{[F]\in\cM(\vv)|\ h^0(F)=2r\}.
\end{equation}
By~(\ref{lowbound}) $U(\vv)$ is open in $\cM(\vv)$.

\begin{lmm}\label{ucompl}{\rm[Markman]}
Keep notation and hypotheses as above. Then $U(\vv)$ is Zariski-dense in
$\cM(\vv)$. Furthermore $U(\vv)=\cM(\vv)$ if
\begin{equation}\label{vsmall}
\langle\vv,\vv\rangle\le(4r-2),
\end{equation}
otherwise
\begin{equation}\label{complcod}
\cod\left(\cM(\vv)\sm U(\vv),\cM(\vv)\right)=2r+1.
\end{equation}
\end{lmm}
\begin{proof}
 Follows from
 Corollary~(3.16), p.~672 of~\cite{markman}. (Notice: the
 definition of $\mu(v)$ is on p.~628, loc.cit.)
 \end{proof}

Let $U^{\flat}(\vv)\ss U(\vv)$ be the open subset defined by
\begin{equation}\label{ubemolle}
U^{\flat}(\vv):=\{[F]\in U(\vv)|\ \mbox{$F$
locally-free and globally generated}\}.
\end{equation}
Let $[F]\in U^{\flat}(\vv)$; the sheaf $E$
appearing in Exact Sequence~(\ref{edef}) is
locally-free, slope-stable and
$v(E)=v(F)^{\vee}$. Thus $[E^{\vee}]\in\cM(\vv)$
and we have a regular map
\begin{equation}\label{prephi}
\begin{matrix}
 U^{\flat}(\vv) & \lra & \cM(\vv)\\
[F] & \mapsto & E^{\vee}
 \end{matrix}
 \end{equation}
\begin{thm}\label{markthm}{\rm[Markman]}
Keep notation and hypotheses as above. There exists an anti-symplectic
(see~(\ref{meno})) birational involution
$$\phi_{\vv}\cl\cM(\vv)\cdots>\cM(\vv)$$
with the following properties:
\begin{itemize}
\item[(1)]
$\phi_{\vv}$ is regular on $U(\vv)$. In particular if~(\ref{vsmall}) holds then
$\phi_{\vv}$ is a regular involution.
\item[(2)]
The restriction of $\phi_{\vv}$ to $U^{\flat}(\vv)$ coincides with the map given
by~(\ref{prephi}).
\item[(3)]
$\phi_{\vv}(U^{\flat}(\vv))=U^{\flat}(\vv)$.
\end{itemize}
\end{thm}
\begin{proof}
 In the notation of~\cite{markman} the map $\phi_{\vv}$ is $\wt{q}_0$ of
 Theorem~(3.21), p.~681, with $a=b=r$ and $\cL$ the line-bundle such that
 $c_1(\cL)=\ell$. Markman does not prove that $\phi_{\vv}$ is anti-symplectic. If
  $r\ge 2$ this follows from Proposition~(\ref{h2action}) below. If $r=1$ the
  map $\phi_{\vv}$ is Beauville's involution and hence it is anti-symplectic
  by Proposition~(\ref{beauvact}). (1) is Item~(1) of Theorem~(3.21), in~\cite{markman}
  (the case $t=0$). (2) is
in~\cite{markman}, first line of p.~683. To prove~(3) it suffices to show that
\begin{equation}\label{dentro}
\phi_{\vv}(U^{\flat}(\vv))\ss U^{\flat}(\vv)
\end{equation}
because $\phi_{\vv}^{-1}=\phi_{\vv}$. Let $[F]\in
U^{\flat}(\vv)$ and let $E$ be the sheaf
appearing in~(\ref{edef}): we must show that
$[E^{\vee}]\in U^{\flat}(\vv)$. We know by
Item~(2) of Lemma~(\ref{stablem}) that $E^{\vee}$
is locally-free. Applying the
$Hom(,\cO_S)$-functor to~(\ref{edef}) we get a
sequence
\begin{equation}\label{dualseq}
0\to F^{\vee}\to H^0(F)^{\vee}\ot\cO_S\to
E^{\vee}\to 0
\end{equation}
which is exact because $F=\wt{F}$ is locally-free. Thus $E^{\vee}$ is
globally-generated. Since $[F]\in U(\vv)$ we have
$H^1(F^{\vee})=H^1(F)^{\vee}=0$ and hence the long exact cohomology sequence
associated to~(\ref{dualseq}) gives $h^0(E^{\vee})=2r$. This
proves~(\ref{dentro}).
 \end{proof}

\begin{lmm}\label{uflatbig}
Keep notation and hypotheses as above and assume
furthermore that $r\ge 2$. Then
$U^{\flat}(\vv)$ is Zariski-dense in $\cM(\vv)$.
\end{lmm}

\begin{proof} Let $\Delta(\vv)\ss\cM(\vv)$ be
given by
\begin{equation}\label{deldef}
\Delta(\vv):=\{[F]\in\cM(\vv)|\ F^{\vee\vee}\not=F\}
\end{equation}
i.e.~the (closed) subset parametrizing singular (not locally-free) sheaves. Let
\begin{equation}\label{tetdef}
\T^0(\vv):=\{[F]\in U(\vv)|\ \wt{F}\not=F\}
\end{equation}
i.e.~the (closed) subset in $U(\vv)$ parametrizing sheaves
 which are \ul{not}
globally generated. By Exact Sequence~(104), p.~683 in~\cite{markman}
\begin{equation}\label{phitet}
\phi_{\vv}(\T^0(\vv))=\Delta(\vv)\cap U(\vv).
\end{equation}
Since $r\ge 2$ we know that $\Delta(\vv)$ is
a proper subset of $\cM(\vv)$ (see
also~(\ref{delteta})). Thus $\T^0(\vv)$ is a proper
subset of $\cM(\vv)$. Since
$U^{\flat}(\vv)=U(\vv)\sm\T^0(\vv)\sm\Delta(\vv)$ we are
done.
 \end{proof}
Lemma~(\ref{ucompl}) and Theorem~(\ref{markthm}) allow us to produce many moduli
spaces $\cM(\vv)$ with a regular anti-symplectic involution. An example: let
$(S,D)$ be a degree-$(2g-2)$ polarized $K3$ and set
$$\vv:=r+c_1(D)+r\eta,\quad g\le r^2+2r.$$
Choosing $D$ as the ample divisor defining (semi)stability of sheaves we see
that the hypotheses of Lemma~(\ref{ucompl}) and Theorem~(\ref{markthm}) are
satisfied except possibly Hypothesis~(\ref{picint}). For $(S,D)$ contained in an
open dense subset of the moduli space of degree-$(2g-2)$ polarized $K3$'s
Hypothesis~(\ref{picint}) is satisfied
 as well
and hence $\phi_{\vv}$ is a regular involution of $\cM(\vv)$. Notice that we get
examples in any (even) dimension.
\subsubsection{Description of $H^2(\phi_{\vv})$ and applications}
Throughout this subsubsection we assume that $\vv$ is given
by~(\ref{vsymmetric}) with $r\ge 2$ and that both Hypothesis~(\ref{picint}) and
~(\ref{eldpos}) hold. By Remark~(\ref{itisgeneric}) we know that $\cM(\vv)$ is a
deformation of $(K3)^{[n]}$ where $2n=2+\la\vv,\vv\ra$, and furthermore
Theorem~(\ref{markthm}) gives us the birational involution
$\phi_{\vv}\cl\cM(\vv)\cdots>\cM(\vv)$.  We also assume that
$\langle\vv,\vv\rangle\ge 2$ i.e.~that $\dim\cM(\vv)\ge 4$. By
Remark~(\ref{itisgeneric}) we know that $\theta_{\vv}\cl\vv^{\bot}\to
H^2(\cM(\vv))$ is an isomorphism of lattices (and Hodge structures). Since
$(\eta-1)\in\vv^{\bot}$ it makes sense to set $h_{\vv}:=\theta_{\vv}(\eta-1)$;
notice that $(h_{\vv},h_{\vv})=2$.  The following result extends to higher rank
the formula of Proposition~(\ref{beauvact}) (notice that if $\vv$ is given
by~(\ref{hilbvect}) then $\theta_{\vv}(\eta-1)=h_g$ where $h_g$ is as
in~(\ref{hgdef})).
\begin{prp}\label{h2action}
Keep notation and hypotheses as above. Then
$H^2(\phi_{\vv})$ is the reflection in the span
of $h_{\vv}$, i.e.
$$H^2(\phi_{\vv})(\a)=-\a+(\a,h_{\vv}) h_{\vv}.$$
\end{prp}
Before proving the proposition we give a
corollary. Let $H_{\vv}$ be the divisor class
such that $c_1(H_{\vv})=h_{\vv}$.

\begin{crl}\label{amplelb}
Keep notation and hypotheses as above, and suppose furthermore that
 $\langle\vv,\vv\rangle\le(4r-2)$, i.e.~that $\dim\cM(\vv)\le 4r$.
Then:
\begin{itemize}
\item[(1)]
$\phi_{\vv}$ extends to all small deformations of
$\cM(\vv)$ that keep $h_{\vv}$ of type $(1,1)$,
\item[(2)]
$H_{\vv}$ is ample.
\end{itemize}
\end{crl}
\begin{proof}
(1): The map $\phi_{\vv}$ is regular by Item~(1)
of Theorem~(\ref{markthm}) and thus Item~(1)
follows from Proposition~(\ref{h2action}) and
Corollary~(\ref{stabinv}). (2): Let $H_0$ be an
ample divisor on $\cM({\vv})$ - it exists because
$\cM({\vv})$ is projective. Since $\phi_{\vv}$ is
regular $\phi_{\vv}^{*}H_0$ is ample. Thus
$(H_0+\phi_{\vv}^{*}H_0)$ is an
 ample divisor class invariant for
 $\phi_{\vv}^{*}$. By
 Proposition~(\ref{h2action}) $(H_0+\phi_{\vv}^{*}H_0)$
 is a multiple of $H_{\vv}$; thus either
 $H_{\vv}$ or $(-H_{\vv})$ is ample. Suppose that
 $(-H_{\vv})$ is ample: we will arrive at a
 contradiction. Let
 $\mu_{\vv}\cl H^2(S)\to
 H^2(\cM(\vv))$ be Donaldson's map
 (see~\cite{ogradyvb}) and $L$ be the
 line-bundle on $S$ such that $c_1(L)=\ell$.
 By Hypothesis~(\ref{picint}) and~(\ref{eldpos})
 $L$ is big and nef and hence
 $$\int_{\cM(\vv)}c_1(-H_{\vv})^{2n-1}
 \wedge\mu_{\vv}(\ell)>0$$
 where $2n=\dim\cM(\vv)$. By
 Fujiki's Formula~(\ref{fujformula}) we get that
 \begin{equation}\label{positivo}
  (-H_{\vv},\mu_{\vv}(\ell))>0.
  \end{equation}
  On the other hand the
 second-to-last
 formula on p.~639 of~\cite{ogradyvb} (warning: the map
 $\theta_{\vv}$ in~\cite{ogradyvb} is the opposite of
 $\theta_{\vv}$ of the present paper) gives that
 $$(-H_{\vv},\mu_{\vv}(\ell))=
 -{1}/{r}\int_S\ell\cdot\ell<0.$$
 This contradicts~(\ref{positivo}); thus $(-H_{\vv})$
 is not ample and hence $H_{\vv}$ must be ample.
\end{proof}

Before proving Proposition~(\ref{h2action}) we prove some lemmas. The first
lemma is very similar to Theorem~(2.9) of~\cite{yoshi}.
\begin{lmm}\label{resact}
Let $U^{\flat}(\vv)$ be as in~(\ref{ubemolle}) and $\iota_{\vv}\cl
U^{\flat}(\vv)\hra \cM(\vv)$ be the inclusion.  Then
$$\iota^{*}_{\vv}\circ H^2(\phi_{\vv})=
\iota^{*}_{\vv}\circ R_{h_{\vv}}.$$
\end{lmm}
 \begin{proof}
 Let $\phi_{\vv}^{\flat}\cl U^{\flat}(\vv)\to U^{\flat}(\vv)$
 be the restriction of $\phi_{\vv}$ (see Item~(3) of
 Theorem~(\ref{markthm})) and $\Phi^{\flat}_{\vv}:=\Id_S\tm\phi_{\vv}^{\flat}$.
Let $\cF$ be the restriction to $S\tm U^{\flat}(\vv)$ of a
 quasi-tautological family on $S\tm\cM(\vv)$. Let $R_{(\eta-1)}\cl H^{*}(S)\to
 H^{*}(S)$ be the reflection in the span of $(\eta-1)$ i.e.
\begin{equation}\label{riflessione}
R_{(\eta-1)}(\a):=-\a+\langle\a,\eta-1\rangle(\eta-1).
\end{equation}
$\vv^{\bot}$ is mapped to itself by $R_{(\eta-1)}$ because
$R_{(\eta-1)}(\vv)=-\vv$. Given $\theta_{\vv}(\a)\in H^2(\cM(\vv))$, where
$\a\in\vv^{\bot}$, we have
\begin{align*}
\iota_{\vv}^{*}\circ H^2(\phi_{\vv})(\theta_{\vv}(\a)) & =\frac{1}{\s(\cF)}
\theta_{(\Phi_{\vv}^{\flat})^{*}\cF}(\a)\\
\iota_{\vv}^{*}\circ R_{h_{\vv}}(\theta_{\vv}(\a) )&
=\frac{1}{\s(\cF)}\theta_{\cF}(R_{(\eta-1)}(\a)),
\end{align*}
hence we must prove that
\begin{equation}\label{cassano}
\theta_{(\Phi_{\vv}^{\flat})^{*}\cF}(\a)=\theta_{\cF} (R_{(\eta-1)}(\a)).
\end{equation}
 Let $\rho\cl S\tm U^{\flat}(\vv)\to U^{\flat}(\vv)$ be the
 projection. By definition of
 $U^{\flat}(\vv)$ and Equation~(\ref{h2andchi})
 we have
 \begin{equation}\label{highvanish}
 R^q\rho_{*}\cF=0,\quad q>0,
 \end{equation}
hence $\rho_{*}\cF$ is
 locally-free of rank $2r\s(\cF)$.
By definition of
 $U^{\flat}(\vv)$ the natural map
$\rho^{*}(\rho_{*}\cF)\to\cF$ is surjective. Let
$\cE$ be the sheaf on $S\tm U^{\flat} (\vv)$
fitting into the exact sequence
\begin{equation}\label{evf}
0\to\cE\to\rho^{*}(\rho_{*}\cF)\to\cF\to 0.
\end{equation}
Let $[F]\in U^{\flat}(\vv)$ and
$[G]=\phi_{\vv}([F])$. By definition of
$\phi_{\vv}$ we have $\cE^{\vee}|_{S\tm[F]}\cong
G^{\s(\cF)}$ and hence the quasi-families
$\cE^{\vee}$, $(\Phi^{\flat}_{\vv})^{*}\cF$ of
sheaves on $S\tm U^{\flat} (\vv)$ with Mukai
vector $\vv$ are equivalent. By
Lemma~(\ref{nondipende}) we get that
\begin{equation}\label{pullphi}
\theta_{(\Phi_{\vv}^{\flat})^{*}\cF}(\a)=
\theta_{\cE^{\vee}}(\a),\quad \a\in\vv^{\bot}.
\end{equation}
By definition of $U^{\flat}(\vv)$ the sheaf $\cF$
is locally-free; taking the dual of~(\ref{evf})
we get that
\begin{equation}\label{totti}
\theta_{\cE^{\vee}}(\a)=
\rho_{*}[\rho^{*}ch(\rho_{*}\cF)^{\vee}\pi^{*}((1+\eta)\a^{\vee})]_6-
\rho_{*}[ch(\cF)^{\vee}\pi^{*}((1+\eta)\a^{\vee})]_6
\end{equation}
where $\pi\cl S\tm U^{\flat}(\vv)\to S$ is the
projection. By~(\ref{highvanish}) we have
$\rho_{!}(\cF)=\rho_{*}(\cF)$ hence
Grothendieck-Riemann-Roch gives that
\begin{equation}\label{grr}
c_1(\rho_{*}\cF)=c_1(\rho_{!}(\cF))=
\rho_{*}[ch(\cF)\pi^{*}(1+2\eta)]_6=\theta_{\cF}(1+\eta).
\end{equation}
Thus
\begin{eqnarray}\label{primo} \rho_{*}
[\rho^{*}ch(\rho_{*}\cF)^{\vee}\pi^{*}((1+\eta)\a^{\vee})]_6
&= & -\rho_{*}[\pi^{*}((1+\eta)\a^{\vee})]_4
c_1(\rho_{*}\cF) \nonumber\\
& = &
\theta_{\cF}(\langle\a,1+\eta\rangle(1+\eta)).
\end{eqnarray}
On the other hand
\begin{equation}\label{secondo}
\rho_{*}[ch(\cF)^{\vee}\pi^{*}((1+\eta)\a^{\vee})]_6=-\theta_{\cF}(\a^{\vee}).
\end{equation}
Plugging~(\ref{primo})-(\ref{secondo})
into~(\ref{totti}) we get that
$$\theta_{\cE^{\vee}}(\a)=
\theta_{\cF}(\a^{\vee}+\langle\a,1+\eta\rangle(1+\eta)).$$
By~(\ref{pullphi})
$$\theta_{(\Phi_{\vv}^{\flat})^{*}\cF}(\a)=
\theta_{\cF}(\a^{\vee}+\langle\a,1+\eta\rangle(1+\eta)).$$
Since $R_{(\eta-1)}(\a)=\a^{\vee}+\langle\a,1+\eta\rangle(1+\eta)$ this
proves~(\ref{cassano}).
 \end{proof}
\begin{lmm}\label{delteta}
Let $\Delta(\vv)$  be given by~(\ref{deldef})
and $\T(\vv)$ be the closure of~(\ref{tetdef}).
\begin{itemize}
\item[(1)]
Both $\Delta(\vv)$ and $\T(\vv)$ are irreducible of codimension $(r-1)$.
\item[(2)]
$\Delta(\vv)\not=\Theta(\vv)$.
\end{itemize}
\end{lmm}

\begin{proof}
(1): It is well-known that $\Delta(\vv)$ is irreducible of codimension $(r-1)$:
it follows from the fact that for any $\ww\in H^{*}(S)$ the moduli space
$\cM(\ww)$ is either empty or of the expected dimension and by irreducibility of
the $Quot$-scheme parametrizing length-$q$ quotients of a fixed locally-free
sheaf on $S$ (Theorem~(6.A.1) of~\cite{huy-lehn}). One also gets that the
generic $F$ parametrized by $\Delta(\vv)$ fits into an exact sequence
\begin{equation}\label{singsheaf}
0\to F\to E\overset{g}{\to} \CC_p\to 0,
\end{equation}
where $[E]\in\cM(\vv+\eta)$ is generic and locally-free, $\CC_p$ is the
skyscraper sheaf at an arbitrary $p\in S$ and $g$ is an arbitrary surjection.
From~(\ref{phitet}) we get that also $\T(\vv)$ is irreducible of codimension
$(r-1)$. (2): We must prove that
\begin{equation}\label{globgen}
\mbox{if $[F]\in\Delta(\vv)$ is generic then $F$
is globally generated.}
\end{equation}
Such an $F$ fits into~(\ref{singsheaf}) where
$[E]\in\cM(\vv+\eta)$ is generic and
locally-free. We claim that  $E$ is globally
generated. The proof is analogous to that of
Lemma~(\ref{uflatbig}). For $[E]$ generic
$h^0(E)=2r+1$ by Corollary~(3.16), p.~672
of~\cite{markman}: let $\wt{E}\ss E$ be the
subsehaf of $E$ generated by global sections,
then $E/\wt{E}$ is Artinian by Lemma~(3.5)
of~\cite{markman}. One considers the Mukai map
$$\begin{matrix}
\cM(\vv+\eta) & \cdots> & M(1+\vv)\\
[E] & \mapsto & [G^{\vee}]
\end{matrix}$$
where $G$ is the sheaf fitting into the exact
sequence
$$0\to G\to H^0(E)\ot\cO_S\to \wt{E}\to 0.$$
A parameter count shows that $\wt{E}=E$ for a generic $[E]\in \cM(\vv+\eta)$, we
leave the details to the reader. Let $[E]\in \cM(\vv+\eta)$ be generic, let
$\pi\cl\PP(E^{\vee})\to S$ be the projection and $\xi$ be the tautological line
sub-bundle of $\pi^{*}E^{\vee}$. Since $H^0(E)=H^0(\xi^{\vee})$ we know that the
linear system $|\xi^{\vee}|$ has no base-locus and hence we have a regular map
$$f\cl \PP(E^{\vee})\to\PP\left(H^0(E)^{\vee}\right)
\cong\PP^{2r}.$$
Let $g\in E^{\vee}_p$ be the map appearing in~(\ref{singsheaf}) and let $x=[g]$;
thus $x\in\PP(E^{\vee})$. If $df(x)$ is injective the sheaf $F$ appearing
in~(\ref{singsheaf}) is globally generated. Thus to prove~(\ref{globgen}) it
suffices to show that $\dim\Im(f)=\dim\PP(E^{\vee})=(r+1)$: this follows from
the easily computed formula
$$\int_{\PP(E^{\vee})}c_1(\xi^{\vee})^{r+1}
=\frac{1}{2}(\ell\cdot\ell)+1.$$
 \end{proof}
Let $U^{\sharp}(\vv)\ss U(\vv)$ be given by
$$U^{\sharp}(\vv):=\{[F]\in U(\vv)|
\ \ell(F^{\vee\vee}/F)+\ell(F/\wt{F})\le 1\}.$$
This is an open subset of $\cM(\vv)$ because both
$$[F]\mapsto\ell(F^{\vee\vee}/F)\quad\mbox{and}\quad
[F]\mapsto\ell(F/\wt{F})$$
are upper semicontinuous functions on the open $U(\vv)$.
\begin{lmm}\label{codsharp}
Keep notation and assumptions as above. Then
$$\cod \left(\cM(\vv)\sm U^{\sharp}(\vv),\cM(\vv)\right)\ge 2.$$
\end{lmm}
\begin{proof}
By Lemma~(\ref{ucompl}) we know that $\cod(\cM(\vv)\sm U(\vv),\cM(\vv))\ge 2$
and hence it suffices to show that $\cod \left(U(\vv)\sm
U^{\sharp}(\vv),\cM(\vv)\right)\ge 2$. We have a decomposition $(U(\vv)\sm
U^{\sharp}(\vv))=A_{(2,0)}\cup A_{(1,1)}\cup A_{(0,2)}$ where
$$\begin{array}{lcl}
A_{(2,0)}& := &\{[F]\in U(\vv)|\ \ell(F^{\vee\vee}/F)\ge 2\},\\
A_{(1,1)}& :=& \{[F]\in U(\vv)|\ \ell(F^{\vee\vee}/F)\ge 1,\ \ell(F/\wt{F})\ge 1\},\\
A_{(0,2)} &:=& \{[F]\in U(\vv)|\ \ell(F/\wt{F})\ge 2\}.
\end{array}$$
$A_{(2,0)}$ is a proper subset of $\Delta(\vv)$, see~(\ref{singsheaf}), and
hence $A_{(2,0)}$ has codimension at least $2$ by Lemma~(\ref{delteta}).
$A_{(1,1)}\ss \Delta(\vv)\cap\Theta(\vv)$ and hence it has codimension at least
$2$ by Lemma~(\ref{delteta}). Finally by~(108) on p.~683 of~\cite{markman} we
have $\phi_{\vv}(A_{(0,2)})\ss A_{(2,0)}$ and hence $A_{(0,2)}$ has codimension
at least $2$.
\end{proof}
 \n
 {\bf Proof of Proposition~(\ref{h2action}).}
 First we prove the proposition in the case $r\ge
 3$. Since
 $U^{\flat}(\vv)=U(\vv)\sm\T(\vv)\sm\Delta(\vv)$
 we get from Item~(1) of
 Lemma~(\ref{delteta}) and Lemma~(\ref{ucompl})
 that
$$\cod(\cM(\vv)\sm U^{\flat}(\vv),\cM(\vv))\ge
r-1.$$
Hence if $r\ge 3$ the map $H^2(\iota_{\vv})$ is
an isomorphism, and thus
Proposition~(\ref{h2action}) follows from
Lemma~(\ref{resact}). We are left with the case
$r=2$, i.e.
$$\vv=2+\ell+2\eta.$$
Let $j_{\vv}\cl U^{\sharp}(\vv)\hra\cM(\vv)$ be
the inclusion. The restriction of $\phi_{\vv}$ to
$U^{\sharp}(\vv)$ is an involution
$\phi_{\vv}^{\sharp}$ of $U^{\sharp}(\vv)$: this
 follows from~\cite{markman}, p.~683. Let
 $\Phi^{\sharp}_{\vv}:=\Id_S\tm\phi_{\vv}^{\sharp}$.
Let $\cF$ be the restriction to $S\tm
U^{\sharp}(\vv)$ of a quasi-tautological family
on $S\tm\cM(\vv)$. Then
\begin{equation}\label{restsharp}
j_{\vv}^{*}\circ H^2(\phi_{\vv})(\theta_{\vv}(\a))=
\theta_{(\Phi_{\vv}^{\sharp})^{*}\cF}(\a), \quad\a\in\vv^{\bot}.
\end{equation}
Let's construct a quasi-family equivalent to $(\Phi_{\vv}^{\sharp})^{*}\cF$. Let
$\wt{\cF}$ be the sheaf on $S\tm U^{\sharp}(\vv)$ given by
$$\wt{\cF}:=\Im\left(\rho^{*}(\rho_{*}\cF)\to\cF\right),$$
where $\rho\cl S\tm U^{\sharp}(\vv)\to U^{\sharp}(\vv)$ is the projection. Let
$\cE$ be the sheaf on $S\tm U^{\sharp}(\vv)$ fitting into the exact sequence
\begin{equation}\label{eftilde}
0\to\cE\to\rho^{*}(\rho_{*}\cF)\to\wt{\cF}\to 0.
\end{equation}
As is easily checked $\cE$ is a quasi-family of torsion-free sheaves on $S$
parametrized by $U^{\sharp}(\vv)$ with Mukai vector $\vv^{\vee}$. If $[F]\in
U^{\sharp}(\vv)$ then $\cE|_{S\tm[F]}\cong (E')^{\s(\cF)}$ where the double dual
of $E'$ is isomorphic to the sheaf $E$ of~(\ref{edef}). $E'$ is slope-stable by
Item~(2) of Lemma~(\ref{stablem}), and furthermore if
$[F]\notin(\T(\vv)\cup\Delta(\vv))$ we have $[E'\ot L]=\phi_{\vv}([F])$, where
$L$ is the line-bundle on $S$ such that $c_1(L)=\ell$: these facts imply that
the quasi-families of sheaves on $S$ with Mukai vector $\vv$ given by
$(\Phi_{\vv}^{\sharp})^{*}\cF$ and $\cE\ot\pi^{*}L$ are equivalent ($\pi\cl S\tm
U^{\sharp}(\vv)\to S$ is the projection). Let $\cG:=\cE\ot\pi^{*}L$.
By~(\ref{restsharp}) and Lemma~(\ref{nondipende}) we have
\begin{equation}\label{gmap}
j_{\vv}^{*}\circ H^2(\phi_{\vv})(\theta_{\vv}(\a))= \theta_{\cG}(\a),
\quad\a\in\vv^{\bot}.
\end{equation}
We compute the right-hand side. We have an exact sequence
\begin{equation}\label{lamseq}
0\to\wt{\cF}\to\cF\to \lambda\to 0
\end{equation}
where $\l$ is a sheaf such that $supp(\l)$ is mapped by $\rho$ to $\T(\vv)\cap
U^{\sharp}(\vv)$ with finite fibers. By Lemma~(\ref{delteta}) $\T(\vv)$ is
irreducible and hence
\begin{equation}\label{pushlambda}
\rho_{*}ch_3(\l)=j_{\vv}^{*}(n c_1(\T(\vv))),\quad n>0.
\end{equation}
Let $\b\in\vv^{\bot}\cap H^{*}(S;\ZZ)$ be the
class such that
\begin{equation}\label{betadef}
\theta_{\vv}(\b)=n c_1(\T(\vv)),
\end{equation}
and let $T_{\b}\cl H^{*}(S)\to H^{*}(S)$ be
defined by
$$T_{\b}(\a):=\a_0\b- \langle e^{-\ell}\wedge\a,
1+\eta \rangle(1+\eta)-e^{-\ell}\wedge\a.$$
Using~(\ref{eftilde}), (\ref{lamseq}),
(\ref{pushlambda}) and~(\ref{grr}) one gets that
$$\theta_{\cG}(\a)=
\theta_{\cF}(T_{\b}(\a)).$$
Now notice that $T_{\b}(\vv^{\bot})=\vv^{\bot}$
(warning: $T_{\b}$ is not an isometry!) and hence
we can rewrite the above equation as
\begin{equation}\label{tetag}
\theta_{\cG}(\a)=j_{\vv}^{*}\circ\theta_{\vv}(T_{\b}(\a)).
\end{equation}
By~(\ref{codsharp}) the map $H^2(j_{\vv})$ is an
isomorphism, hence~(\ref{tetag}) together
with~(\ref{gmap}) gives that
\begin{equation}\label{rosibindi}
 H^2(\phi_{\vv})(\theta_{\vv}(\a))=\theta_{\vv}(T_{\b}(\a)),\quad
 \a\in\vv^{\bot}.
\end{equation}
Since~(\ref{utile}) is a homomorphism $H^2(\phi_{\vv})$ is an isometric
involution and hence the restriction of $T_{\b}$ to $\vv^{\bot}$ is an isometric
involution. We claim that this implies that
\begin{equation}\label{betexp}
\b=-(g-3)-\ell-2\eta.
\end{equation}
Since with this value of $\b$ we have
$T_{\b}=R_{(\eta-1)}$ the above equation proves
Proposition~(\ref{h2action}) for $r=2$. Let's
prove that~(\ref{betexp}) holds. Let
$\L_{\vv}:=(\CC\op\CC\ell\op\CC\eta)\cap\vv^{\bot}$
and $\Xi_{\vv}:=H^2(S)\cap\ell^{\bot}$; thus
$\vv^{\bot}=\L_{\vv}\op_{\bot}\Xi_{\vv}$. Since
$T_{\b}(\Xi_{\vv})=\Xi_{\vv}$ we have
$T_{\b}(\L_{\vv})=\L_{\vv}$; since
$(1-\eta)\in\L_{\vv}$ and
$T_{\b}(1-\eta)\in\b+\L_{\vv}$ we have
$\b\in\L_{\vv}$. Let $\b=b_0+b_2\ell+b_4\eta$
where $b_i\in\ZZ$: thus
\begin{equation}\label{mastella}
 T_{\b}(\a)= \a_0
b_0+\a_0 b_2\ell+\a_0 b_4\eta - \langle
e^{-\ell}\wedge\a, 1+\eta
\rangle(1+\eta)-e^{-\ell}\wedge\a.
\end{equation}
A straightforward computation shows that the restriction of the above map to
$\vv^{\bot}$ is an involution only if~(\ref{betexp}) holds.
 \qed

\msk
 It follows from~(\ref{betadef}) and~(\ref{betexp}) that
\begin{equation}\label{teteq}
 c_1(\T(\vv)) =  \theta_{\vv}(-(g-3)-\ell-2\eta)
 \qquad\mbox{if $r=2$.}
\end{equation}
By~(\ref{phitet}) we have $\Delta(\vv)=\phi_{\vv}^{*}\T(\vv)$ and hence
\begin{equation}\label{deleq}
 c_1(\Delta(\vv))  =
\theta_{\vv}(2+\ell+(g-3)\eta)\qquad\mbox{if $r=2$.}
\end{equation}
\subsection{Another example}\label{dualexample}
 \setcounter{equation}{0}
Let $S$ be a  $K3$ surface, $H_S$ be a divisor on $S$ with
\begin{equation}\label{hons}
H_S\cdot H_S=2g-2.
\end{equation}
Let $h_S:=c_1(H_S)$. Let $\mu$ and $\xi_2$ be as in~(\ref{mumap})
and~(\ref{h2hilb}) respectively and let $h\in H^{1,1}_{\ZZ}(S^{[2]})$ be
\begin{equation}\label{testh}
h:=\mu(h_S)-t\xi_2.
\end{equation}
 If
\begin{equation}\label{gcond}
g=2+t^2
\end{equation}
then $(h,h)=2$ and hence the I Conjecture~(\ref{iconj}) predicts the existence
of an anti-symplectic birational involution of $S^{[2]}$. Let us test the
conjecture for $t=0,1,2$. It suffices to produce the anti-symplectic involution
under the hypothesis that $(S,H_S)$ is the generic couple with $S$ a $K3$ and
$H_S$ an ample degree-$(2g-2)$ divisor of a given divisibility. In fact once
this is proved the degeneration procedure of Section~(\ref{limpliesi}) will give
that the desired involution exists in general. Assume that $(S,H_S)$ is
polarized and generic; thus by~\cite{mayer} the generic curve in $|H_S|$ is
smooth and its genus is given by the $g$ appearing in~(\ref{hons}). If $t=0$
then $g=2$ hence $S$ is a double cover of $\PP^2$ and the covering involution
$\phi_S\cl S\to S$ induces an anti-symplectic involution $\phi$ of $S^{[2]}$.
Notice that  $H^2(\phi)$ is never equal to $R_h$ (the reflection in $\ZZ h$)
because $\phi^{*}\Delta_2=\Delta_2$ where $\Delta_2$ is as in~(\ref{nonred}). If
$t=1$ then $g=3$. Let $\phi\cl S^{[2]}\to S^{[2]}$ be Beauville's involution
defined in Subsection~(\ref{beauvex}). By Proposition~(\ref{beauvact}) we have
$H^2(\phi)=R_h$ generically. (See also the comments at the end of
Subsection~(\ref{beauvex}).) Now let $t=2$. In this case we get a new example.
Since $H_S$ is ample and $g=6$ we know by~\cite{mayer} that the linear system
$|H_S|$ has no base-locus and that it defines an embedding
$S\hra|H_S|^{\vee}\cong\PP^6$ as a degree-$10$ surface. According to
Mukai~\cite{mukaik3} the generic $K3$ surface of degree $10$ in $\PP^6$ is
described as follows. Let $\GG (1,\PP^4)\ss\PP^9$ be the Pl\"ucker embedding of
the Grassmannian of lines in $\PP^4$ and let
\begin{equation}\label{fano3fold}
F:=\GG (1,\PP^4)\cap\PP^6
 \end{equation}
 where
$\PP^6\ss\PP^9$ is a linear space transversal to $\GG (1,\PP^4)$. We have
$K_F\cong\cO_F(-2)$ and $\deg F=5$ i.e.~$F$ is a Fano $3$-fold of index $2$ and
degree $5$; by Iskovskih(\cite{iskovskih}, Thm.~(4.2)) there is one  isomorphism
class of such Fano $3$-folds. If $\ov{Q}\ss\PP^6$ is a quadric transversal to
$F$ then
\begin{equation}\label{sfq}
S:=F\cap \ov{Q}
 \end{equation}
 is a $K3$ surface of degree $10$ in $\PP^6$. Mukai~(\cite{mukaik3}, Cor.~(4.3))
proved that the generic $K3$ of degree $10$ is given by~(\ref{sfq}). We will
define an involution on $S^{[2]}$ analogous to Beauville's involution on
$T^{[2]}$ where $T\ss\PP^3$ is a quartic surface: one replaces $\PP^3$ by $F$
and $\GG(1,\PP^3)$ by the Hilbert scheme $W(F)$ parametrizing conics in $F$.
Before defining the involution we state two results on lines and conics lying on
$F$. Let $R(F)$ be the Hilbert scheme parametrizing lines contained in $F$. For
$p\in F$ let
\begin{equation}\label{}
R_p:=\{[L]\in R(F)|\ p\in L\}.
\end{equation}
The following result is due to Iskovskih.
\begin{lmm}[Cor.~(6.6) of \cite{iskovskih}]\label{rnice}
Keeping notation as above, $R(F)$ is isomorphic to $\PP^2$.
\end{lmm}
For $Z\in(\PP^6)^{[2]}$ we let $\la Z\ra\ss\PP^6$ be the line spanned by $Z$.
Let $B\ss F^{[2]}$ be the closed subset given by
\begin{equation}\label{bdef}
B:=\{[Z]\in F^{[2]}|\ \la Z\ra\ss F\}.
\end{equation}
For $[Z]\in F^{[2]}$ and $[\ov{L}]\in R(F)$ let
\begin{align*}
W_{Z}:=&\{[C]\in W(F)|\ Z\ss C\},\nonumber\\
R_{\ov{L}}:=&\{[L]\in R(F)|\ [\ov{L}+L]\in W(F)\}.
\end{align*}
We will prove the following.
\begin{lmm}\label{wnice}
Keep notation as above. Then
\begin{itemize}
\item[(1)]
Let $[\ov{L}]\in R(F)$. Then $R_{\ov{L}}$ is isomorphic to $\PP^1$.
\item[(2)]
Let $[Z]\in F^{[2]}$. If $[Z]\in B$ then $W_Z$
consists of those $[C]\in W(F)$ such that $C$
contains $\langle Z\rangle$ (hence $W_Z$ is
identified with $R_{\la Z\ra}$). If $[Z]\notin B$
then $W_{Z}$ consists of a single conic.
\item[(3)]
$W(F)$ is irreducible, rational of dimension $4$.
\end{itemize}
\end{lmm}
Granting this lemma for the moment being, we
proceed to define the involution. We make the
following assumption regarding $S$:
\begin{equation}\label{nolnoc}
\mbox{$S$ contains no line and no conic.}
\end{equation}
Let $B_S:=B\cap S^{[2]}$ - this makes sense because we have an inclusion
$S^{[2]}\ss F^{[2]}$. Let $U:=\left(S^{[2]}\sm B_S\right)$;  $U$ is
Zariski-dense in $S^{[2]}$ because $\cod(B_S,S^{[2]})=2$
 by Lemma~(\ref{rnice}). We will define a regular map
 \begin{equation}\label{phireg}
\phi_U\cl U\to S^{[2]}.
 \end{equation}
Let $[Z]\in U$. Since $[Z]\notin B_S$ there is a unique conic $C_Z\ss F$
containing $Z$ (Item~(2) of Lemma~(\ref{wnice})). As is easily checked the ideal
of $Z$ in $\cO_{C_Z}$ is locally principal because $[Z]\notin B_S$. Thus there
is a well-defined residual scheme $Z'$ of $Z\ss (C_Z\cap S)$ in $C_Z$. The
intersection $C_Z\cap S$ is a scheme of length $4$ by~(\ref{nolnoc}), and hence
$[Z']\in S^{[2]}$.  We define the map $\phi_U$ of~(\ref{phireg}) by setting
$\phi_U([Z]):=([Z'])$.
\begin{prp}\label{geomact}
Keep notation as above. Then:
\begin{itemize}
\item[(1)]
 $\phi_U$ is regular and it extends to a birational involution $\phi\cl
S^{[2]}\cdots>S^{[2]}$,
\item[(2)]
$H^2(\phi)=R_h$ where $h=\mu(h_S)-2\xi_2$.
\end{itemize}
\end{prp}
\begin{proof}
(1): One verifies easily that $\phi_U$ is
regular. Furthermore $\phi_U(U)=U$ and $\phi_U$
is an  involution of $U$: thus $\phi$ is a
birational involution.  (2): By~(\ref{sfq}) and
the hypothesis the complement of $U$ in $S^{[2]}$
has codimension $2$: since $\phi$ is regular on
$U$ it follows that $H^2(\phi)$ is locally
constant over the family of $K3$ surfaces given
by~(\ref{sfq}). Hence arguing as in the proof of
Proposition~(\ref{beauvact}) we see  that  it
suffices to show that:
\begin{itemize}
\item[(2a)]
$\phi^{*}h=h$,
\item[(2b)]
$\phi^{*}\Delta_2\not=\Delta_2$,
\item[(2c)]
$\phi^{*}\s^{[2]}=-\s^{[2]}$,
\end{itemize}
where $\s^{[2]}$ is the symplectic form induced
on $S^{[2]}$ by a symplectic form $\s$ on $S$.
(2a): Consider the regular map
\begin{equation}\label{dualf}
\begin{matrix}
S^{[2]} & \overset{f}{\lra} &
|I_S(2)|^{\vee}\cong\PP^5\\
[Z] & \mapsto & \{Q\in|I_S(2)|\ |\ \langle Z\rangle\ss Q\}.
\end{matrix}
\end{equation}
 One checks easily that
$f^{*}c_1\left(\cO_{\PP^5}(1)\right)=h$; since $f$ commutes with
$\phi$ this proves Item~(2a). (2b): If the quadric $\ov{Q}$
of~(\ref{sfq}) is chosen generically then the generic conic
parametrized by $W(F)$ which is tangent to $\ov{Q}$ intersects
$\ov{Q}$ in two other distinct points and hence
$\phi^{*}\Delta_2\not=\Delta_2$. (2c): Let $[Z_i]\in U$ and
$\phi([Z_i])=[Z'_i]$ for $i=1,2$: by Lemma~(\ref{wnice}) $W(F)$ is
rational and hence
$$Z_1+Z'_1\sim Z_2+Z'_2$$
where $\sim$ is rational equivalence.
Thus $\phi^{*}\s^{[2]}=-\s^{[2]}$ by
Mumford's Theorem on $0$-cycles
(see Prop~(22.24) of~\cite{voibook}).
\end{proof}
\n
 {\bf Proof of Lemma~(\ref{wnice}).}
We need to recall Iskovskih's description of a birational map between $F$ and a
smooth quadric $3$-fold ((6.5), p.511 of~\cite{iskovskih}). There exists a line
$L_0\ss F$ such that $N_{L_0/F}\cong \cO_{L_0}\op\cO_{L_0}$. Projection from
$L_0$ defines a birational map
$$\pi\cl F\cdots> Q_0$$
to a smooth quadric $3$-fold $Q_0$. The indeterminacy locus of $\pi$ is resolved
by blowing up $L_0$. There is a smooth twisted cubic $Y\ss Q_0$ such that the
induced regular map $Bl_{L_0}(F)\to Q_0$ contracts the proper transform of any
line in $F$ meeting $L_0$ to a point of $Y$. The inverse $\pi^{-1}\cl Q_0\cdots>
F$ is given by $|I_Y(2)|$, the linear system of quadrics containing $Y$. The
indeterminacy locus of $\pi^{-1}$ is resolved by blowing up $Y$ and the regular
map $Bl_Y(Q_0)\to F$ contracts the proper transform of every chord of $Y$
contained in $Q_0$ (the intersection $\langle Y\rangle\cap Q_0$ is a smooth
quadric surface). In order to prove~(1) of Lemma~(\ref{wnice}) we describe
$R(F)$ via $\pi$. Let $[L]\in (R(F)\sm\{[L_0]\})$. From the above it follows
that if $L\cap L_0=\es$ then $\pi(L)\ss Q_0$ is a line meeting $Y$ and not
contained in $\langle Y\rangle$, and that viceversa every such line is mapped by
$\pi^{-1}$ to a line $L\ss F$ with $L\cap L_0=\es$. If $L\cap L_0\not=\es$ then
$\pi$ contracts $L$ to a point and viceversa to every point of $Y$ there
corresponds an $[L]\in R(F)$. (Associating to a point $q\in Y$ the unique line
$R\ss\langle Y\rangle\cap Q_0$ such that $R\cap Y=q$ we may view the latter case
as a degeneration of the former.) Item~(1) of Lemma~(\ref{wnice}) follows easily
from the above description of $R(F)$. In order to prove~(2) we describe $W(F)$
via $\pi$. Let $[C]\in W(F)$. If $C\cap L_0=\es$ then $\pi(C)$ is a conic not
contained in $\langle Y\rangle$ with $\pi(C)\cap Y$ a scheme of length $2$, and
viceversa every conic in $Q_0$ not contained in $\langle Y\rangle$ and
intersecting $Y$ in a subscheme of length $2$ is mapped by $\pi^{-1}$ to a conic
$C\ss F$ with $C\cap L_0=\es$. If $C\cap L_0$ is a point then $\pi$ maps $C$ to
a line in $Q_0$ not contained in $\langle Y\rangle$, which meets $Y$ if and only
if $C$ is reducible, and viceversa every line in $Q_0$ and not contained in
$\langle Y\rangle$ is mapped by $\pi^{-1}$ to a conic $C\ss F$ such that $C\cap
L_0$ is a point. (The latter case can be viewed as a degeneration of the former
case by adding to the line in $Q_0$ a suitable chord of $Y$.) Of course if
$C\supset L_0$ then $f(C)$ is a point of $Y$ and viceversa.... Item~(2)  of
Lemma~(\ref{wnice}) follows easily from the above description of $W(F)$. Let's
prove Item~(3). By the above description $W(F)$ is irreducible of dimension $4$
and the generic conic parametrized by $W(F)$ is smooth. Let $H\ss F$ be a smooth
hyperplane section: thus $H$ is a del Pezzo surface of degree $5$ - the blow up
of $\PP^2$ at $4$ points no $3$ of which are collinear. Since the set of conics
belonging to $H$ has dimension $1$ we have a rational map
$$\begin{matrix}
W(F) & \cdots> & H^{[2]}\\
[C] & \mapsto & [C\cap H].
\end{matrix}$$
Let $L_1,\ldots,L_{10}\ss H$ be the $10$ lines of
$H$: by Item~(2) of Lemma~(\ref{wnice}) the above
map is an isomorphism on the complement of
$\bigcup_{i=1}^{10}L_i^{[2]}$. Thus $W(F)$ is
birational to $H^{[2]}$ and hence it is rational.
 \qed
\subsection{Two involutions}\label{dueinv}
 \setcounter{equation}{0}
Let $X$ be an irreducible symplectic variety with $H_1,H_2$ ample
divisors. Let $h_i:=c_1(H_i)$; we assume that $h_1,h_2$ are
linearly independent. Suppose that there exist regular involutions
$\phi_i\cl X\to X$ for $i=1,2$ such that $H^2(\phi_i)=R_{h_i}$.
Let $\psi:=\phi_1\circ\phi_2$; then $H^2(\psi)$ is described as
follows. Let $\L:=\RR h_1\op\RR h_2$. The restriction of
Beauville's form to $\L$ has signature $(1,1)$ and hence we have a
direct sum decomposition
$$H^2(X;\RR)=\L\op\L^{\bot}.$$
There exist a basis $\{e_{+},e_{-}\}$ of $\L$ and a real number
$\l>1$ such that
$$(e_{+},e_{+})=(e_{-},e_{-})=0,\quad H^2(\psi)(e_{\pm})=\l^{\pm
1}.$$
In particular $\psi$ has infinite order and it should give an
interesting dynamical system on $X$ (see the Introduction
of~\cite{mcmullen}). Now assume that $X$ and $\phi_1,\phi_2$ are
defined over a number filed $K$. Following Silverman~\cite{silver}
(now we should assume, and this is always possible, that $h_i\in
(\RR_{+}e_{+}\op \RR_{+}e_{-})$)  we may associate to $e_{\pm}$ a
normalized logarithmic height $h_{\pm}$ defined on all $p\in
X(\ov{K})$ and having the following properties:
\begin{itemize}
\item[(1)]
 $h_{\pm}(p)\ge 0$ for all $p\in X(\ov{K})$,
\item[(2)]
 $h_{\pm}(\psi(p))=\l^{\pm 1}h_{\pm}(p)$,
\item[(3)]
 $h_{+}(p)=h_{-}(p)=0$ if and only if the orbit
 $\{\psi^j(p)\}_{j\in\ZZ}$ is finite.
\end{itemize}
This might be used to produce many rational points on $X$. Silverman's argument
shows that there are no $\psi$-invariant effective non-zero cycles of dimension
or codimension $1$. Thus if $p\in X(\ov{K})$ with $h_{+}(p)\not=0$ or
$h_{-}(p)\not=0$ then the Zariski-closure of $\{\psi^j(p)\}_{j\in\ZZ}$ is all of
$X$ or some effective cycle of dimension $1<d<(\dim X-1)$. A concrete example:
Let $S\ss\PP^3\tm\PP^3$ be the complete intersection of $\Si_1,\ldots,\Si_4\in
|\cO_{\PP^3}(1)\boxtimes\cO_{\PP^3}(1)|$. For general $\Si_1,\ldots,\Si_4$ the
surface $S$ is a $K3$  and the projection  $\pi_i\cl S\to\PP^3$ to the $i$-th
factor is an isomorphism to a smooth  quartic with no lines. Let
$\ell_i:=c_1(\pi_i^{*}\cO_{\PP^3}(1))$. Let $X:=S^{[2]}$ and
$h_i:=(\mu(\ell_i)-\xi_2)$ where $\mu$ and $\xi_2$ are as in~(\ref{h2hilb});
this is an example of the situation described above. In fact we have Beauville's
involution $\phi_i\cl S^{[2]}\to S^{[2]}$ with $H^2(\phi_i)=R_{h_i}$, see
Subsection~(\ref{beauvex}). Assume that $\Si_1,\ldots,\Si_4$ are defined over a
number field $K$: we can show that there exists a finite extension $K'\supset K$
with $X(K')$ Zariski-dense proceeding as follows. There exists a $K'$ such that
we have a curve $C\in|\pi_1^{*}\cO_{\PP^3}(1)|$ defined over $K'$ which is
birational to $\PP^1_{K'}$. Then
$$R:=\{[Z]\in S^{[2]}|\ supp(Z)\in C,\quad [Z]\in\Delta_2\}$$
is a ruled surface with $R(K')$ Zariski-dense in it. One checks that the
$\psi$-orbit of the Poincar\'e dual of $R$ is infinite. Since there are no
$\psi$-invariant effective non-zero divisors we get that the Zariski
 closure of $\{\psi^j(R(K'))\}_{j\in\ZZ}$ is the whole of $S^{[2]}$.
\section{Examples: linear systems}\label{exlinsyst}
\setcounter{equation}{0}
We will give examples
 of degree-$2$ polarized
deformations of $(K3)^{[n]}$ which satisfy the hypotheses of
Proposition~(\ref{stablin}), i.e.~evidence in favor of the L
Conjecture~(\ref{lconj}). The examples are inspired by Mukai
(Ex.(5.17)~\cite{mukaisug}): he gave the $4$-dimensional example. In proving
that the linear systems of our examples are well behaved we will verify that a
so-called {\it Strange duality} (see~\cite{dontu,lepot,danila,lepot2}) holds for
the linear systems in question. We will make the connection with Strange duality
in a separate subsection. In the last subsection we will examine more closely
the $4$-dimensional example.
\subsection{The examples}\label{guisa}
\setcounter{equation}{0}
Let $S$ be a $K3$ surface, $D$ an ample divisor on $S$ and let
\begin{equation}\label{vranktwo}
\vv:=2+\ell+2\eta,\quad \la\vv,\vv\ra\le 6.
\end{equation}
We assume that both Hypothesis~(\ref{picint}) and Inequality~(\ref{eldpos}) hold
- recall that this is always possible, see the end of
Subsubsection~(\ref{mukairefl}). Thus by Item~(1) of Theorem~(\ref{markthm}) the
involution $\phi_{\vv}\cl\cM(\vv)\to\cM(\vv)$ is regular and by
Corollary~(\ref{amplelb}) the divisor class $H_{\vv}$ is ample. We will show
that Proposition~(\ref{stablin}) holds for $X_0=\cM(\vv)$ and $H_0=H_{\vv}$. Let
$L$ be the line-bundle on $S$ such that $c_1(L)=\ell$: we assume that
\begin{equation}\label{lample}
\mbox{$L$ is ample,}\quad L\cdot L=2g-2.
\end{equation}
By~(\ref{vranktwo}) we have
\begin{equation}\label{dimcong}
\dim\cM(\vv)=2(g-4),\quad    g\le 8.
\end{equation}
It follows from our hypotheses and the results of
Mayer~\cite{mayer} that $|L|$ has no base-locus
and that the map
$$f_S\cl S\to |L|^{\vee}\cong\PP^g$$
is an embedding. From now on $S$ is embedded in $\PP^g$ by the map $f_S$. If
$g=8$ we make the following extra assumption. Let
$${\mathbb Gr}(2,\CC^6)\hra \PP(\wedge^2\CC^6)\cong\PP^{14}$$
be the Pl\"ucker embedding: we assume that
\begin{equation}\label{voilas}
S={\mathbb Gr}(2,\CC^6)\cap\PP^8,\qquad\mbox{$\PP^8$ transversal to ${\mathbb
Gr}(2,\CC^6)$}.
\end{equation}
The generic polarized $K3$ of degree $14$ is obtained in this way,
see~\cite{mukaik3}. (This assumption forces us to choose $L,D$ with $L^{\ot k}
\cong\cO_S(D)$.) Let $\Si\ss|I_S(2)|$ be the locally closed subset given by
$$\Si:=\{Q\in|I_S(2)|\,\,\,|\quad \rk Q=6,\quad sing(Q)\cap
S=\es\}.$$
If $Q\in\Si$ let $F_{g-3}(Q)\ss{\mathbb Gr}(g-3,\PP^g)$ be the subset
parametrizing $(g-3)$-dimensional linear subspaces contained in $Q$; thus
$F_{g-3}(Q)$ has two connected components. If $\L\in F_{g-3}(Q)$ the
intersection number $(\L\cdot S)_Q$ of $\L,S$ as cycles on $Q$ is well-defined
because $S$ is contained in the smooth locus of $Q$. If $\L,\L'\in F_{g-3}(Q)$
belong to the same connected component then $(\L\cdot S)_Q=(\L'\cdot S)_Q$. On
the other hand if $\L,\L'$ belong to different connected components then
$$(\L\cdot S)_Q+(\L'\cdot S)_Q=2g-2.$$
Thus if $Q\in \Si$ there exists an integer $0\le i(Q)\le g-1$ such that for
every $\L\in F_{g-3}(Q)$ we have
$$(\L\cdot S)_Q=(g-1)\pm i(Q).$$
Let
$$\Si_a:=\{Q\in\Si|\quad i(Q)=a\}.$$
Each $\Si_a$ is an open subset of $\Si$ and we have
$\Si=\Si_0\cup\cdots\cup\Si_{g-1}$ (disjoint union). Let $Y:=\ov{\Si}_0$.
 Now we are ready to describe $f_{\vv}\cl\cM(\vv)\cdots>|H_{\vv}|^{\vee}$.
\begin{prp}\label{explmap}
Keep notation and assumptions as above.  Then:
\begin{itemize}
\item[(1)]
$|H_{\vv}|$ has no base-locus.
\item[(2)]
There is an isomorphism $|H_{\vv}|^{\vee}\cong |I_S(2)|$ such that
$\Im(f_{\vv})=Y$.
\item[(3)]
The map $f_{\vv}\cl \cM(\vv)\to Y$ is finite of degree $2$; the corresponding
covering involution is equal to $\phi_{\vv}$.
\end{itemize}
\end{prp}
Since $\left(\cM(\vv),H_{\vv}\right)$ is a degree-$2$ polarized deformation of
$(K3)^{[n]}$ the above proposition gives examples in favour of the L
Conjecture~(\ref{lconj}) for $\dim=4,6,8$. In fact by Corollary~(\ref{stablin})
we get that in each of these dimensions the L Conjecture holds for at least one
irreducible component of the relevant moduli space. The proof of the above
proposition actually identifies $\cM(\vv)$ with a natural \lq\lq double
cover\rq\rq of $Y$. We explain this.
\begin{clm}\label{atleastfive}
If $Q\in|I_S(2)|$ then $\rk(Q)\ge 5$.
\end{clm}
\begin{proof}
If $\rk(Q)\le 4$ there exists a reducible hyperplane section of $S$,
contradicting Hypothesis~(\ref{picint}).
\end{proof}
Let $Q\in Y$: by the above claim $\rk(Q)=5$ or $\rk(Q)=6$ and hence $F_{g-3}(Q)$
has one or two connected components respectively. ($F_{g-3}(Q)$ is defined as
above also if $\rk(Q)=5$.) Let $\cF_{g-3}\to Y$ be the map with fiber
$F_{g-3}(Q)$ over $Q$ and
$$\cF_{g-3}\lra W\overset{\zeta}{\lra}Y,$$
be its Stein factorization, thus $\zeta^{-1}(Q)$ is the set of connected
components of $F_{g-3}(Q)$: then $\zeta\cl W\to Y$ is finite of degree two. Let
$\rho\cl W^{\nu}\to W$ be the normalization map.
\begin{prp}\label{eccemoduli}
Keep notation and assumptions as above. The map $f_{\vv}\cl \cM(\vv)\to Y$ lifts
to a (regular) map $\tilde{f}_{\vv}\cl \cM(\vv)\to W$. The map $\cM(\vv)\to
W^{\nu}$ induced by $\tilde{f}_{\vv}$ is an isomorphism.
\end{prp}
\subsection{Proof of~(5.1)-(5.3)}
\setcounter{equation}{0}
\subsubsection{Sheaves on $S$ and quadrics in $|I_S(2)|$}
Let $F$ be a $D$-slope-stable sheaf on $S$ with
\begin{equation}\label{condfa}
v(F)=2+\ell+s\eta.
\end{equation}
Let $A\ss H^0(F)$ be a subspace. Let
\begin{equation}\label{deltaf}
\e_A\cl \wedge^2 A\to H^0(\wedge ^2 F)\cong H^0(\cO_S(1))
\end{equation}
be the natural map. Now assume that $\dim A\ge 3$. Let $F_A\ss F$ be the sheaf
generated by $A$. By Lemma~(3.5) of~\cite{markman} the quotient $F/F_A$ is
Artinian; let
 $\O_A:=supp(F/F_A)\cup sing(F)$.
 Thus $\O_A$ is a finite set of points.
 Let $E_A$ be the locally-free sheaf on $S$
 fitting into the exact sequence
 \begin{equation}\label{rieccola}
 0\to E_A\to A\ot\cO_S\to F_A\to 0.
 \end{equation}
 We have a regular map
 \begin{equation}\label{mapgr}
\begin{array}{ccl}
 S\sm\O_A &
\overset{\l_A}{\lra} & {\mathbb
Gr}\left(2,A^{\vee}\right)
\ss\PP\left(\wedge^2 A^{\vee}\right)\\
x & \mapsto & Ann(E_A)_x
\end{array}
\end{equation}
Now suppose that $\dim A=4$: since the Grassmannian above is a quadric
hypersurface we get by \lq\lq pull-back\rq\rq a quadric in $|I_S(2)|$. To be
precise  choose a trivialization $\wedge^4 A^{\vee}\overset{\sim}{\to}\CC$ and
let $R_A\in Sym^2\left(\wedge^2 A\right)$ correspond to multiplication on
$\wedge^{2}A^{\vee}$.  Let $P_A:=Sym^2(\e_A)(R_A)$: since $R_A$ vanishes on
${\mathbb Gr}\left(2,A^{\vee}\right)$ we have
$$P_A\in \Ker\left(Sym^2 H^0(\cO_S(1))\to H^0(\cO_S(2))
\right).$$
\begin{lmm}\label{nonzero}
Keep notation and assumptions as above, in particular assume that $\dim A=4$.
Then $P_A\not=0$.
\end{lmm}
\begin{proof}
We claim that
\begin{equation}\label{menodiuno}
\dim(\Ker\e_A)\le 1.
\end{equation}
The lemma follows from the above inequality
because $R_A$ is non-degenerate. Let $\s,\tau\in
A$: we claim that
\begin{equation}\label{decinj}
\mbox{if $\e_A(\s\wedge\tau)=0$ then
$\s\wedge\tau=0$.}
\end{equation}
If $\s=0$ there is nothing to prove so we may
assume that $\s\not=0$. By $D$-slope-stability of
$F$ and Hypothesis~(\ref{picint}) $\s$ defines a
map $\s\cl\cO_S\to F$ which is injective on
fibres away from a finite set $Z\ss S$; thus away
from $Z$ we have $\tau=f\s$ for a regular
function $f$. Since $\cod(Z,S)=2$ the function
$f$ extends to a regular function on $S$ and
hence is equal to a constant $c$; since $F$ is
torsion-free we get that $\tau=c\s$.
Now~(\ref{decinj}) gives that
$\PP(\Ker\e_A)\cap{\mathbb Gr}(2,A)=\es$; since
${\mathbb Gr}(2,A)$ is a hypersurface in
$\PP(\wedge^2 A)$ this implies~(\ref{menodiuno}).
\end{proof}
\begin{dfn}\label{qadef}
Keep notation and assumptions as above and suppose that $\dim A=4$.  Then
$Q_A:=V(P_A)$ is a well-defined quadric hypersurface in $\PP^g$ by
Lemma~(\ref{nonzero}). Clearly $Q_A\in |I_S(2)|$ and $\rk(Q_A)\le 6$. When
$h^0(F)=4$ we set $Q_F:=Q_{H^0(F)}$.
\end{dfn}
\subsubsection{The map $q_{\vv}\cl\cM(\vv)\to Y$ and its
lift to $W$}
Let $[F]\in\cM(\vv)$. By our hypotheses and
 Lemma~(\ref{ucompl}) we know that $h^0(F)=4$.
Thus we have a map
 \begin{equation}\label{qvdef}
 \begin{matrix}
 \cM(\vv) & \overset{q_{\vv}}{\lra} & |I_S(2)|\\
 [F] & \mapsto & Q_F.
 \end{matrix}
 \end{equation}
A moment's thought shows that $q_{\vv}$ is
regular.
 \begin{lmm}\label{immagine}
 Keeping notation as above, $\Im(q_{\vv})= Y$.
 \end{lmm}
\begin{proof}
First we prove that
\begin{equation}\label{contenuto}
\Im(q_{\vv})\ss Y.
\end{equation}
Let $U^{\flat}(\vv)\ss\cM(\vv)$ be the open
subset defined in~(\ref{ubemolle}). Let
\begin{equation}\label{condsuf}
[F]\in
 U^{\flat}(\vv),\quad \phi_{\vv}([F])\not=[F].
 \end{equation}
 We will show that $Q_F\in \Si_0$.
This will prove~(\ref{contenuto}): in fact $U^{\flat}(\vv)$ is dense in
$\cM(\vv)$ by Lemma~(\ref{uflatbig}), and $\phi_{\vv}$ is not the identity (it
is anti-symplectic!), hence the set of $[F]$ satisfying~(\ref{condsuf}) is
Zariski-dense in $\cM(\vv)$. Let $\G:=H^0(F)$; since $[F]\in U^{\flat}(\vv)$ we
have $F_{\G}=F$. By Theorem~(\ref{markthm}) we have
$\phi_{\vv}([F])=[E^{\vee}_{\G}]$. Let's show that $\rk Q_F=6$, i.e.~that
$\Ker(\e_{\G})=0$. Assume the contrary and let $\a\in\Ker(\e_{\G})$ be non-zero.
Let $K$ be the sheaf on $S$ given by $K:=\Ker\left(\wedge^2
\G\ot\cO_S\to\wedge^2 F \right)$: since $\e_{\G}(\a)=0$ we have $\a\in H^0(K)$.
We have a natural exact sequence
$$0\to \wedge^2 E_{\G}\to K\overset{\pi}{\to}E_{\G}\ot
F\to0.$$
By~(\ref{decinj}) we know that $\a$ is a rank-$4$ element of $\wedge^2 \G$ and
hence
$$\pi(\a)\in H^0(E_{\G}\ot F)=\Hom(E_{\G}^{\vee},F)$$
is an isomorphism. Since $\phi_{\vv}([F])=[E^{\vee}_{\G}]$ this
contradicts~(\ref{condsuf}): thus $\rk Q_F=6$. Let's show that $sing(Q_F)\cap
S=\es$. The restriction to $S\ss\PP^g$ of projection from $sing(Q_F)$ is
identified with $\l_{\G}$. If $sing(Q_F)\cap S\not=\es$ then $sing(Q_F)\cap S$
is $0$-dimensional by Hypothesis~(\ref{picint}) and hence $\l_{\G}$ is not
regular; since $[F]\in U^{\flat}(\vv)$ we have $\O_{\G}=\es$ hence $\l_{\G}$ is
defined on all of $S$, contradiction. Thus $Q_F\in\Si$. Let's show that
$Q_F\in\Si_0$. Choose a non-zero $\s\in H^0(F)$ and let
$$\G_{\s}:=
\{V\in {\mathbb Gr}\left(2,H^0(F)^{\vee}\right)|\
f(\s)=0,\ \forall f\in V\}$$
be the corresponding Schubert cycle. There exists a unique $\L_{\s}\in
F_{g-3}(Q_F)$ which gives $\G_{\s}$ when projected from $sing(Q_F)$. We have
$\left(\L_{\s}\cdot S\right)_{Q_F}=\deg(\s)$ where $(\s)$ is the zero-locus of
$\s$ - notice that $(\s)$ is $0$-dimensional by Hypothesis~(\ref{picint}). It
follows from~(\ref{vranktwo})-(\ref{lample}) that $c_2(F)=(g-1)\eta$ and hence
$\deg(\s)=(g-1)$. This shows that $Q_F\in\Si_0$ and thus
proves~(\ref{contenuto}). Now let's prove that $Y\ss\Im(q_{\vv})$. Let $Q\in
\Si^0$: we will show that there exists $[F]\in\cM(\vv)$ such that $Q_F=Q$.
Projecting $S$ from $sing(Q)\cong\PP^{g-6}$ we get a morphism $\l\cl S\to\PP^5$
with $\l(S)\ss\ov{Q}$, where $\ov{Q}\ss\PP^5$ is the image of $Q$, a smooth
quadric hypersurface. Choose an isomorphism $\ov{Q}\cong {\mathbb
Gr}\left(2,\CC^4\right)$ and let $\xi$ be the tautological  bundle on ${\mathbb
Gr}\left(2,\CC^4\right)$. Then $F:=\l^{*}\xi^{\vee}$ is a vector-bundle on $S$
with
$$\rk(F)=2,\quad c_1(F)=\ell,\quad c_2(F)=(g-1)\eta.$$
(The last equality holds because $Q\in\Si_0$.) Since $F$ is globally generated
and since $h^2(F)=0$ (this is easily checked) we get by Lemma~(3.5)
of~\cite{markman} ((3)$\implies$(2))  that $F$ is $D$-slope-stable. Thus
$[F]\in\cM(\vv)$. Since $S\ss\PP^g$ is non-degenerate so is $S\ss\PP^5$: it
follows that $Q_F=Q$. This proves that $\Si_0\ss\Im(q_{\vv})$: since
$\Im(q_{\vv})$ is closed we get that $Y\ss\Im(q_{\vv})$.
\end{proof}
The map $q_{\vv}$ lifts to a map $\tilde{q}_{\vv}\cl\cM(\vv)\to W$ almost by
construction. Given $[F]\in\cM(\vv)$ we can associate to $F$ not only the
quadric $Q_F$ but also a choice of component of $F_{g-3}(Q_F)$: if $\rk(Q_F)=6$
then $F_{g-3}(Q_F)$ is naturally isomorphic to the variety parametrizing planes
in ${\mathbb Gr}\left(2,H^0(F)^{\vee}\right)$ and hence it is clear how to
choose a component of $F_{g-3}(Q_F)$, if $\rk(Q_F)=5$ there is only one
component to choose. As is easily checked $\tilde{q}_{\vv}$ is regular. Let
$\tilde{q}^{\nu}_{\vv}\cl\cM(\vv)\to W^{\nu}$ be the (regular) map induced by
$\tilde{q}_{\vv}$ (recall that $\rho\cl W^{\nu}\to W$ is the normalization map).
\begin{lmm}\label{uguali}
Keep notation as above. Then $\tilde{q}^{\nu}_{\vv}$ is an isomorphism, and
hence $q_{\vv}\cl\cM(\vv)\to Y$ has degree $2$. Via $\tilde{q}^{\nu}_{\vv}$ the
 involution of $W^{\nu}$ defined by the degree-$2$ finite map $W^{\nu}\to Y$
 coincides with $\phi_{\vv}$.
\end{lmm}
\begin{proof}
We prove the first statement. Clearly $\tilde{q}^{\nu}_{\vv}$ is an isomorphism
over $\rho^{-1}(\Si_0)$.  Hence it suffices to verify that
$\tilde{q}^{\nu}_{\vv}$ has finite fibers . Let $H$ be a hyperplane section of
the projective space $|I_S(2)|$: since $q_{\vv}\circ\phi_{\vv}=q_{\vv}$ we have
$\phi_{\vv}^{*}(q_{\vv}^{*}H)=q_{\vv}^{*}H$. Thus by
Proposition~(\ref{h2action}) we have $q_{\vv}^{*}H=k H_{\vv}$ for some
$k\in\ZZ$. Since $\dim Y\not=0$ (in fact $\dim Y=\dim\cM(\vv)$) we get
$k\not=0$, and since $H_{\vv}$ is ample by Corollary~(\ref{amplelb}) we have
\begin{equation}\label{multiplo}
q_{\vv}^{*}H\sim k H_{\vv}, \quad k>0.
\end{equation}
If $\tilde{q}^{\nu}_{\vv}$ has a positive dimensional fiber then $q_{\vv}^{*}H$
is trivial on such a fiber: by the above equality also $H_{\vv}$ is trivial on
that fiber, contradicting ampleness of $H_{\vv}$. Let's prove the second
statement of the lemma. Let $\iota_{\vv}$ be the
 involution of $W^{\nu}$ defined by the degree-$2$ finite map $W^{\nu}\to Y$:
  over $\rho^{-1}(\Si_0)$ we clearly have $\iota_{\vv}=
\tilde{q}^{\nu}_{\vv}\circ\phi_{\vv}\circ(\tilde{q}^{\nu}_{\vv})^{-1}$, and
hence the same equality holds on all of $W^{\nu}$.
\end{proof}
\subsubsection{An isomorphism
$|H_{\vv}|\overset{\sim}{=} |I_S(2)|^{\vee}$}
\begin{lmm}\label{kappa1}
Keep notation and assumptions as above. Then $q_{\vv}^{*}H\sim H_{\vv}$.
\end{lmm}
\begin{proof}
If $\dim\cM(\vv)=0$, i.e.~$g=4$, there is nothing to prove. If
$\dim\cM(\vv)=2$, i.e.~$g=5$, the result is immediate. In fact in
this case $Y=|I_S(2)|\cong\PP^2$ and hence $\deg(q_{\vv}^{*}H\cdot
q_{\vv}^{*}H)=\deg(q_{\vv})=2$; since $\deg(H_{\vv}\cdot
H_{\vv})=2$ the result follows from~(\ref{multiplo}). Assume that
$\dim\cM(\vv)\ge 4$, i.e.~that $g\ge 6$. We will prove the lemma
by computing the intersection numbers $\deg(q_{\vv}^{*}H\cdot R)$
and $\deg(H_{\vv}\cdot R)$ where $R$ is a certain curve in
$\cM(\vv)$ parametrizing singular sheaves. First we construct the
curve. Let
\begin{equation}\label{wdef}
\ww:=2+\ell+3\eta.
 \end{equation}
Since $g\ge 6$ we have $\cM(\ww)\not=\es$ by Theorem~(\ref{yoshi}). There exists
$[V]\in\cM(\ww)$ such that
\begin{equation}\label{sceltav}
\mbox{$V$ is locally-free globally generated, $h^0(V)=\chi(V)=5$,}
 \end{equation}
 see the proof of Lemma~(\ref{delteta}). Choose $p\in S$ and let
$R:=\PP(V^{\vee}_p)$. Let $\pi_S\cl S\tm R\to S$ be the projection and $\iota\cl
R\hra S\tm R$ be the inclusion $\iota(x):=(p,x)$. Let $\cF$ be the sheaf on
$S\tm R$ fitting into the exact sequence
\begin{equation}\label{dustin}
0\to\cF\lra \pi_S^{*}V\overset{\a}{\lra}
\iota_{*}\cO_R(1)\to 0,
\end{equation}
where $\a$ is induced by the tautological quotient
$V_p\ot\cO_R\to\cO_R(1)$. For
$x\in R$ let $F_x:=\cF|_{S\tm\{x\}}$. Then $F_x$ is torsion-free, $v(F_x)=\vv$,
and $F_x$ is $D$-slope-stable because $V$ is $D$-slope-stable. Since $F_x\cong
F_{x'}$ only if $x=x'$ the $R$-flat family $\cF$ induces an inclusion
$R\hra\cM(\vv)$. Let's prove that
\begin{equation}\label{hvr}
\deg(H_{\vv}\cdot R)=1.
\end{equation}
Since $c_1(H_{\vv})=\theta_{\vv}(\eta-1)$ we have
\begin{equation}\label{cicikov}
\deg(H_{\vv}\cdot
R)=\deg\theta_{\cF}(\eta-1)=-\int_{S\tm R}
ch_3(\cF)=\int_{S\tm R}ch_3(\iota_{*}\cO_R(1)).
\end{equation}
Applying Grothendieck-Riemann-Roch one gets that
the last term equals $1$; this
proves~(\ref{hvr}). Now let's prove that
\begin{equation}\label{natavota}
\deg(q_{\vv}^{*}H\cdot R)=1.
\end{equation}
We must describe the restriction of $q_{\vv}$ to $R$. Let
\begin{equation}\label{betamap}
\begin{matrix}
\PP\left(H^0(V)^{\vee}\right) & \overset{\b_V}{\hra} & |I_S(2)| \\
x & \mapsto & Q_{Ann(x)}
\end{matrix}
\end{equation}
where $Q_{Ann(x)}$ is as in Definition~(\ref{qadef}) - this makes sense because
since $h^0(V)=5$ we have $\dim Ann(x)=4$. Let $\G:=H^0(V)$ and let $\l_{\G}$ be
as in~(\ref{mapgr}). Then $V_p^{\vee}$ is naturally identified with $\l_{\G}(p)$
and hence
$$R=\PP(V_p^{\vee})=\PP(\l_{\G}(p))\ss\PP(H^0(V)^{\vee}).$$
Let $x\in R$: it is immediate that
\begin{equation}\label{qvbv}
q_{\vv}([F_x])=\b_V(x).
 \end{equation}
 Since $\b_V$ is linear
Equation~(\ref{natavota}) follows. The lemma follows
from~(\ref{multiplo}),
(\ref{hvr}) and~(\ref{natavota}).
\end{proof}
Since $S$ is projectively normal
Hirzebruch-Riemann-Roch gives that
\begin{equation}\label{granderoma}
\dim|I_S(2)|= d(g):=\frac{1}{2}(g-2)(g-3)-1.
\end{equation}
\begin{lmm}\label{nondegenere}
Keep notation and assumptions as above. Then $Y$ is a non-degenerate subvariety
of $|I_S(2)|$.
\end{lmm}
\begin{proof}
If $\dim\cM(\vv)=0$, i.e.~$g=4$, then $|I_S(2)|$ is a point and
the result is trivially true. If $\dim\cM(\vv)=2$ i.e.~$g=5$ then
$Y=|I_S(2)|$ and again the result holds. If $\dim\cM(\vv)=4$
i.e.~$g=6$ then $\dim|I_S(2)|=5$ and of course $\dim Y=4$. Thus if
$Y$ is degenerate it is a hyperplane, and since $\deg q_{\vv}=2$
we get by Lemma~(\ref{kappa1}) that $\int_{\cM(\vv)}h^4_{\vv}=2$.
On the other hand Fujiki's formula~(\ref{fujformula}) and the
value of Fujiki's constant for deformations of $(K3)^{[2]}$
(see~(\ref{fujconst})) give that $\int_{\cM(\vv)}h^4_{\vv}=12$,
contradiction. Now assume that $\dim\cM(\vv)=6$ i.e.~$g=7$. Let
$\ww$ be as in~(\ref{wdef}) and $[V]\in\cM(\ww)$
satisfying~(\ref{sceltav}). Let $\G:=H^0(V)$ and
\begin{equation}\label{tvdef}
\Delta_V:=\bigcup_{p\in
S}\PP(\l_{\G}(p))\ss\PP\left(H^0(V)^{\vee}\right),
\end{equation}
where $\l_{\G}$ is as in~(\ref{mapgr}); thus $\Delta_V$ is the
image of $\PP(V^{\vee})$ under the natural map
$\zeta_V\cl\PP(V^{\vee})\to\PP\left(H^0(V)^{\vee} \right)$.
Proceeding exactly as in the proof of Lemma~(\ref{kappa1}) we can
associate to each $x\in\PP(V^{\vee})$  a $D$-slope-stable
torsion-free singular sheaf $F_x$ on $S$ with $v(F_x)=\vv$: this
gives an inclusion $\PP(V^{\vee})\ss\cM(\vv)$. By~(\ref{qvbv}) we
have
\begin{equation}\label{cupiello}
q_{\vv}|_{\PP(V^{\vee})}=\b_V\circ\zeta_V.
\end{equation}
Since  $\Delta_V$ spans $\PP\left(H^0(V)^{\vee} \right)$ and since
$\b_V$ is linear we get that
$$\la Y\ra=\la \Im(q_{\vv})\ra \supset\Im(\b_V)$$
where $\la Y\ra$ is the span of $Y$. Choose $[V']\in\cM(\ww)$ with
$[V']\not=[V]$ and $V'$ locally-free globally generated with
$h^0(V')=\chi(V')=5$. Thus
\begin{equation}\label{contiene}
\la Y\ra \supset\Im(\b_V)\cup\Im(\b_V').
\end{equation}
We claim that $\Im(\b_V)$ and $\Im(\b_{V'})$ are disjoint linear spaces. The
lemma for $g=7$ follows from this because in this case $\dim|I_S(2)|=9$
by~(\ref{granderoma}) and on the other hand from~(\ref{contiene}) we get that
$\dim\la Y\ra\ge 9$. Let's prove that $\Im(\b_V)$ and $\Im(\b_{V'})$ are
disjoint. Assume that there exists a quadric $Q\in\Im(\b_V)\cap \Im(\b_{V'})$.
Projection in $\PP^g$ with center $sing(Q)$ defines a rational map $\rho\cl
S\cdots>\PP^k$ where $k=5$ if $\rk(Q)=6$ and $k=4$ if $\rk(Q)=5$. In the former
case $\Im(\rho)$ is contained in a smooth quadric hypersurface which we identify
with ${\mathbb Gr}(2,\CC^4)$: let $\xi$ be the tautological vector-bundle on
${\mathbb Gr}(2,\CC^4)$. Let $J:=S\cap sing(Q)$; thus $J$ is $0$-dimensional by
Hypothesis~(\ref{picint}). By definition of $\b_V,\b_{V'}$ the restriction  of
$\xi^{\vee}$ to $(S\sm J)$ is isomorphic to $V|_{S\sm J}$ and to $V'|_{S\sm J}$.
Thus $V|_{S\sm J}\cong V'|_{S\sm J}$: since $\cod(J,S)=2$ and since $V,V'$ are
locally-free the isomorphism extends to all of $S$ and hence $[V]=[V']$,
contradiction. If $\rk(Q)=5$ then $\Im(\rho)\ss\PP^4$ is contained in a smooth
quadric hypersurface $\ov{Q}$; embedding $\ov{Q}$ in ${\mathbb Gr}(2,\CC^4)$ as
a hyperplane section for the Pl\"ucker embedding we may proceed as in the
previous case and we will again arrive at a contradiction. This finishes the
proof of the lemma when $g=7$. Finally assume that $\dim\cM(\vv)=8$ i.e.~$g=8$.
We will show that the subset of $\cM(\vv)$ parametrizing sheaves $F$ with
$F^{\vee\vee}/F$ of length $2$ has image in $Y$ spanning all of $|I_S(2)|$. Let
$\uu:=2+\ell+4\eta$. Then $\la\uu,\uu\ra=-2$ hence by Theorem~(\ref{yoshi}) the
moduli space $\cM(\uu)$ consists of a single point $[W]$. Arguing as in the
proof of Lemma~(\ref{delteta}) one shows that $W$ is locally-free globally
generated and that $h^0(W)=\chi(W)=6$. Choose a trivialization of $\wedge^6
H^0(W)^{\vee}$ and let
$$\begin{matrix}
 \wedge^2 H^0(W)^{\vee} & \lra & Sym^2
(\wedge^2 H^0(W))\\
y & \mapsto & R_y
 \end{matrix}$$
 be defined by $R_y(\a,\b):=
(y\wedge\a\wedge\b)$ for $\a,\b\in \wedge^2 H^0(W)^{\vee}$. Let ${\mathbb
Gr}:={\mathbb Gr}(2,H^0(W)^{\vee})$ and $\cO_{\mathbb Gr}(1)$ be the Pl\"ucker
line-bundle. Then
$$R_y\in\Ker
\left(Sym^2 H^0(\cO_{\mathbb Gr}(1))\to
H^0(\cO_{\mathbb Gr}(2)\right).$$
Let $\G:=H^0(W)$ and let $\e_{\G}$ be as
in~(\ref{deltaf}); then
$$P_y:=Sym^2(\e_{\G})(R_y)\in
\Ker\left(Sym^2 H^0(\cO_S(1))\to H^0(\cO_S(2)\right).$$
\begin{clm}\label{annamo}
Keep notation and hypotheses as above. If $P_y=0$ then $y=0$.
\end{clm}
\begin{proof}
Recall that $S$ is given by~(\ref{voilas}). Let $\xi$ be the tautological
rank-two vector-bundle on ${\mathbb Gr}(2,\CC^6)$. The vector-bundle
$\xi^{\vee}|_S$ is globally generated and as is easily checked
$h^2(\xi^{\vee}|_S)=0$, thus by Lemma~(3.5) of~\cite{markman} ((3)$\implies$(2))
we get that $\xi^{\vee}|_S$ is $D$-slope-stable. An easy computation gives
$v(\xi^{\vee}|_S)=\uu$ and thus $\xi^{\vee}|_S\cong W$. Hence
$V(P_y)=V(R_y)\cap\PP^8$. Since $\PP^8$ is transversal to ${\mathbb Gr}$ the
restriction map
$$H^0(I_{\mathbb Gr}(2))\to H^0(I_S(2))$$
is injective: this proves the claim.
\end{proof}
Since $\dim|\cO_{\mathbb Gr}(1)|^{\vee}=14=\dim |I_S(2)|$
(use~(\ref{granderoma})) we get an isomorphism
\begin{equation}\label{deltaiso}
\begin{array}{rrcc}
\d\cl & |\cO_{\mathbb Gr}(1)|^{\vee}=\PP(\wedge^2 H^0(W)^{\vee}) & \overset{\sim}{\lra} & |I_S(2)|\\
 & [y]\quad\quad\ \  & \mapsto & [P_y]
\end{array}
\end{equation}
Let $\Delta_W\ss\PP(H^0(W)^{\vee})$ be defined as
in~(\ref{tvdef}). Let $x,x'\in \Delta_W$ be distinct points, thus
$x\in\PP(\l_{\G}(p))=\PP(W_p^{\vee})$,
$x'\in\PP(\l_{\G}(p'))=\PP(W_{p'}^{\vee})$.
 Let $F_{x,x'}$ be the sheaf on $S$ fitting into the exact sequence
$$0\to F_{x,x'}\to V\overset{\pi}{\to}\CC_p\op\CC_{p'}\to 0,$$
where $\pi$ is determined by $x,x'$. Then $F_{x,x'}$ is torsion-free
$D$-slope-stable and $v(F_{x,x'})=\vv$; thus $[F_{x,x'}]\in\cM(\vv)$. Similarly
to~(\ref{qvbv}) we have that
$$q_{\vv}([F_{x,x'}])=\d(\la x,x'\ra).$$
Here $\la x,x'\ra\in{\mathbb Gr}\ss|\cO_{\mathbb Gr}(1)|^{\vee}$. Thus
\begin{equation}\label{corda}
\la Y\ra=\la\Im(q_{\vv})\ra\supset \la \d(chord(\Delta_W))\ra,
\end{equation}
where $chord(\Delta_W)$ is the set of chords of $\Delta_W$,
naturally embedded in ${\mathbb Gr}\ss|\cO_{\mathbb
Gr}(1)|^{\vee}$. Since $\Delta_W$ is a non-degenerate $3$-fold in
$\PP\left(H^0(W)^{\vee} \right)\cong\PP^5$ it follows  that
$chord(\Delta_W)$ is non-degenerate in $|\cO_{\mathbb
Gr}(1)|^{\vee}$. Since $\d$ is an isomorphism this proves the
lemma when $g=8$.
\end{proof}
\begin{prp}\label{identificazione}
Keep notation and assumptions as above. The map $q_{\vv}\cl\cM(\vv)\to|I_S(2)|$
defines an isomorphism $q_{\vv}^{*}\cl |I_S(2)|^{\vee}\overset{\sim}{\to}
|H_{\vv}|$.
\end{prp}
\begin{proof}
By Lemmas~(\ref{kappa1})-(\ref{nondegenere}) we have an injective linear map
$q_{\vv}^{*}\cl |I_S(2)|^{\vee}\hra |H_{\vv}|$. By Lemma~(\ref{amplelb}) the
divisor $H_{\vv}$ is ample; since $K_{\cM(\vv)}\sim 0$ we have
$\dim|H_{\vv}|=\left(\chi(\cO_{\cM(\vv)}(H_{\vv}))-1\right)$. Let
$\la\vv,\vv\ra+2=2n$; since $\cM(\vv)$ is a deformation of $(K3)^{[n]}$ one
knows, see p.96 of~\cite{elgole}, that
\begin{equation}\label{hkrr}
\chi(\cO_{\cM(\vv)}(H_{\vv}))={n+2\choose 2}.
\end{equation}
Substituting $n=(g-4)$ and using~(\ref{granderoma}) we get that $\dim
|I_S(2)|^{\vee}=\dim |H_{\vv}|$, hence $q_{\vv}^{*}$ is an isomorphism.
\end{proof}
\subsubsection{Proof of
Propositions~(\ref{explmap})-(\ref{eccemoduli}).}
Proposition~(\ref{identificazione}) identifies $f_{\vv}\cl\cM(\vv)\to
|H_{\vv}|^{\vee}$ with
 $q_{\vv}\cl\cM(\vv)\to |I_S(2)|$. Since $q_{\vv}$ has no base points and
 $\Im(q_{\vv})=Y$ we get Items~(1)-(2) of Proposition~(\ref{explmap}).
 Lemma~(\ref{uguali}) gives Item~(3) of Proposition~(\ref{explmap}) and also
 Proposition~(\ref{eccemoduli})
\qed
\subsection{Strange duality}
\setcounter{equation}{0}
The isomorphism
\begin{equation}\label{ourstrdua}
|H_{\vv}|^{\vee}\cong |I_S(2)|
\end{equation}
of Item~(2) of Proposition~(\ref{explmap}) can be interpreted as a particular
case of a conjectural Strange duality  between spaces of sections of determinant
line-bundles on moduli
 spaces of sheaves on curves~\cite{dontu} or surfaces~\cite{lepot,danila,lepot2}. We
 we will formulate the strange duality statement and then we will make the connection
 with~(\ref{ourstrdua}). Let $(S,D)$ be a polarized $K3$. For $i=0,1$ let
 $$\vv_i:=r_i+\ell_i+s_i\eta\in H^0(S;\ZZ)_{\ge 0}\op H^{1,1}_{\bf Z}(S)\op
 H^4(S;\ZZ)$$
and let $\cM(\vv_i)$ be the moduli space of pure sheaves $F$ on $S$ with
$v(F)=\vv_i$ and semi-stable with respect to $D$. We assume that $D$ is
$\vv_i$-generic for $i=0,1$ and that $(r_i+\ell_i)$ is indivisible, so that
Theorems~(\ref{yoshi})-(\ref{h2iso}) apply. Assume furthermore that
 \begin{equation}\label{vperp}
\la \vv_i,\vv_{1-i}^{\vee}\ra  =0,
 \end{equation}
and that
\begin{align}
(r_0\ell_1+r_1\ell_0)\cdot D & > 0,\label{segnopos}\\
\mbox{or $(r_0\ell_1+r_1\ell_0)\cdot D$} & < 0,\label{segnoneg}
\end{align}
Let $\cL(\vv_{1-i})$ be the (holomorphic) line-bundle on $\cM(\vv_i)$ such that
$$c_1(\cL(\vv_{1-i}))=\begin{cases}
\theta_{\vv_i}(-\vv_{1-i}^{\vee}) & \mbox{if~(\ref{segnopos}) holds},\\
\theta_{\vv_i}(\vv_{1-i}^{\vee}) & \mbox{if~(\ref{segnoneg}) holds.}
\end{cases}$$
Mimicking the proof of Theorem~(2.1) of~\cite{danila} one gets the following.
\begin{prp}
Keep notation and assumptions as above. There is a section $\s_{\vv_1,\vv_0}\in
H^0\left(\cM(\vv_0)\tm\cM(\vv_1); \cL(\vv_1)\boxtimes\cL(\vv_0)\right)$,
canonical up to multiplication by a non-zero scalar, such that
$$(\s_{\vv_1,\vv_0})=\{([E_0],[E_1])\in \cM(\vv_0)\tm\cM(\vv_1)|\ h^1(E_0\ot E_1)>0\}.$$
\end{prp}
\begin{proof}
We assume that there exists a tautological sheaf $\cF(\vv_i)$ on
$S\tm\cM(\vv_i)$ for $i=0,1$; if this is not the case one works with the
tautological sheaf on $S\tm Quot$ where $Quot$ is a suitable Quot-scheme and
then applies a descent argument (see~\cite{danila}). Let $\rho_i\cl
S\tm\cM(\vv_0)\tm\cM(\vv_1)\to S\tm\cM(\vv_i)$ and $\pi\cl
S\tm\cM(\vv_0)\tm\cM(\vv_1)\to \cM(\vv_0)\tm\cM(\vv_1)$ be the projections. We
consider the line-bundle $\cL$ on $\cM(\vv_0)\tm\cM(\vv_1)$ defined by
$$\cL:=\det\pi_{!}\left(\rho_0^{*}\cF(\vv_0)\ot\rho_1^{*}\cF(\vv_1)\right).$$
Let $[F_i]\in\cM(\vv_i)$. Applying Grothendieck-Riemann-Roch we get that
$$c_1(\cL|_{[F_0]\tm\cM(\vv_1)})=\theta_{\vv_1}(\vv_0^{\vee}),\qquad
c_1(\cL|_{\cM(\vv_0)\tm [F_1]})=\theta_{\vv_0}(\vv_1^{\vee}).$$
(Recall that $\vv_i^{\vee}\in\vv_{1-i}^{\bot}$ by~(\ref{vperp}).) Since
$\cM(\vv_i)$ has $b_1=0$ it follows that
\begin{equation}\label{eccoelle}
\cL\cong
\begin{cases}
\cL(\vv_1)^{-1}\boxtimes\cL(\vv_0)^{-1} & \mbox{if~(\ref{segnopos}) holds,}\\
\cL(\vv_1)\boxtimes\cL(\vv_0) & \mbox{if~(\ref{segnoneg}) holds.}
\end{cases}
\end{equation}
By~(\ref{vperp}) we have $\chi(F_0\ot F_1)=0$ and furthermore
\begin{align}
H^2(F_0\ot F_1)=0 &\mbox{\ \ if~(\ref{segnopos}) holds,}\\
H^0(F_0\ot F_1)=0 &\mbox{\ \ if~(\ref{segnoneg}) holds.}
\end{align}
It follows by standard arguments that there exists a canonical section
$$\s\in
\begin{cases}
H^0(\cL^{-1}) &\mbox{if~(\ref{segnopos}) holds,}\\
H^0(\cL) &\mbox{if~(\ref{segnoneg}) holds,}
\end{cases}$$
such that $\s([F_0],[F_1])=0$ if and only if $h^1(F_0\ot F_1)>0$.
\end{proof}
We may view $\s_{\vv_1,\vv_0}$ as a map
\begin{equation}\label{strdualmap}
\s_{\vv_1,\vv_0}\cl H^0(\cM(\vv_0);\cL(\vv_1))^{\vee}\to
H^0(\cM(\vv_1);\cL(\vv_0)).
\end{equation}
\begin{stt}\label{strdua}{\rm [Strange duality]}
The map $\s_{\vv_1,\vv_0}$ of~(\ref{strdualmap}) is an isomorphism.
\end{stt}
A comment: by Theorem~(\ref{yoshi}) $\cM(\vv_i)$ is a deformation of
$(K3)^{[n_i]}$ , where $2n_i=2+\la\vv_i,\vv_i\ra$. By a well-known formula, see
p.96 of~\cite{elgole}, we have
$$\chi(\cM(\vv_i);\cL(\vv_{1-i}))=
{\frac{1}{2}(c_1(\cL(\vv_{1-i})),c_1(\cL(\vv_{1-i})))+n_i+1\choose n_i}.$$
By Theorem~(\ref{h2iso}) we know that $\theta_{\vv_i}$ is an isometry and hence
we get that
$$\chi(\cM(\vv_i);\cL(\vv_{1-i}))=
{n_i+n_{1-i}\choose n_i}.$$
Hence if $\cL(\vv_0)$ and $\cL(\vv_1)$ have no higher cohomology
there spaces of global sections have equal dimensions: this is
consistent with the Strange duality statement. Now we show that
the map~(\ref{ourstrdua}) is the projectivization
of~(\ref{strdualmap}) for suitable $\vv_0,\vv_1$. Let
\begin{equation}\label{vedoppiav}
  \vv:=r+\ell+r\eta,\qquad \ww:=1-\eta
\end{equation}
and assume that Hypothesis~(\ref{picint}) and~(\ref{eldpos}) hold. Letting
$\vv_0:=\vv$, $\vv_1:=\ww$ we see that all of our previous assumptions are
verifed; let's spell out the Strange duality statement in this case, in
particular for $r=2$. We have
\begin{equation}
c_1(\cL(\ww))=h_{\vv}:=\theta_{\vv}(\eta-1).
\end{equation}
We let $H_{\vv}$ be a divisor on $\cM(\vv)$ such that $c_1(H_{\vv})=h_{\vv}$. On
the other hand $\cM(\ww)=S^{[2]}$ and
\begin{equation}\label{esplicito}
c_1(\cL(\vv))=\mu(\ell)-r\xi_2
\end{equation}
where $\mu$ and $\xi_2$ are as in Subsubsection~(\ref{hilbdescr}). Now set
$r=2$, let $L$ be the line-bundle such that $c_1(L)=\ell$ and let $L\cdot
L=(2g-2)$. We assume that $L$ is very ample and hence $S\ss\PP^g$ is
non-degenerate with $L\cong\cO_S(1)$. We claim that there is a canonical
identification
\begin{equation}\label{beddo}
|c_1(\cL(\vv))|\cong |I_S(2)|.
\end{equation}
In fact if $\mu(L)$ is the line-bundle on $S^{[2]}$ such that
$c_1(\mu(L))=\mu(\ell)$ we have a canonical identification $|\mu(L)|\cong
|\cO_{\PP^g}(2)|$ and as is easily verified
$|c_1(\cL(\vv))|=|\mu(L)(-\Delta_2)|$ is the subsystem $|I_S(2)|$. Furthermore
$f_{\ww}\cl S^{[2]}\cdots>|I_S(2)|^{\vee}$ is identified with the regular map
\begin{equation}\label{fdoppiav}
 \begin{matrix}
S^{[2]} & \overset{f_{\ww}}{\lra} & |I_S(2)|^{\vee}\\
[Z] & \mapsto & \{Q|\ Q\supset\la Z\ra\}
\end{matrix}
\end{equation}
($S$ contains no lines and is cut out by quadrics because of
Hypothesis~(\ref{picint}), and hence $f_{\ww}$ is regular.) Thus
statement~(\ref{strdua}) asserts that $\s_{\vv,\ww}$ gives an
isomorphism $|I_S(2)|^{\vee}\cong |H_{\vv}|$.
Proposition~(\ref{explmap}) gives such an isomorphism for $g\le 8$
- to be precise when $g=8$ we have the extra
assumption~(\ref{voilas}). Let's show that
Isomorphism~(\ref{ourstrdua}) is equal to the projectivization of
$\s_{\vv,\ww}$. Let $q_{\vv}$ be the map of~(\ref{qvdef}).
\begin{clm}\label{cesemo}
Keep notation as above. Let $[F]\in\cM(\vv)$ and assume that
$q_{\vv}([F])\in\Si_0$. Let $[Z]\in S^{[2]}$. Then
$q_{\vv}([F])\in f_{\ww}([Z])$ if and only if $h^1(I_Z\ot F)>0$.
\end{clm}
\begin{proof}
By definition $q_{\vv}([F])\in f_{\ww}([Z])$ if and only if $Q_F\supset\la
Z\ra$. By the definition of $Q_F$ (see Definition~(\ref{qadef})) this is
equivalent to the existence of a non-zero $\tau\in H^0(F)$ vanishing on $Z$,
i.e.~to $h^0(I_Z\ot F)>0$. Since $\chi(I_Z\ot F)$ and $h^2(I_Z\ot F)$ both
vanish the claim follows.
\end{proof}
Now let $\g\cl |I_S(2)|\overset{\sim}{\to} |H_{\vv}|^{\vee}$ be the inverse of
the isomorphism of Proposition~(\ref{explmap}) and consider the composition
$$\g\circ\s_{\ww,\vv}\cl |H_{\vv}|^{\vee}\cdots> |H_{\vv}|^{\vee},$$
a priori a rational linear map. Let $x\in\Si_0$: by Claim~(\ref{cesemo}) we know
that $\g\circ\s_{\ww,\vv}$ is regular at $x$ and that
$\g\circ\s_{\ww,\vv}(x)=x$. By Lemma~(\ref{nondegenere}) $\Si_0$ is
non-degenerate and hence $\g\circ\s_{\ww,\vv}$ is regular everywhere and equal
to the identity.
\subsection{The $4$-dimensional case }
\setcounter{equation}{0}
Keep notation and assumptions as in the introduction to
Subsection~(\ref{guisa}) and suppose that $g=6$. Then
by~(\ref{dimcong}) we have $\dim\cM(\vv)=4$;  set $X:=\cM(\vv)$.
We will present a couple of observations on $Y=f_{\vv}(X)$. We
assume that $S$ is the generic $K3$ of genus $g$ (this forces
$\cO_S(D)\cong L^{\ot k}$) and hence
\begin{equation}\label{eccoesse}
S=F\cap\ov{Q}
\end{equation}
where $F$ is the Fano $3$-fold given by~(\ref{fano3fold}) and
$\ov{Q}$ is a quadric hypersurface transversal to $F$,
see~(\ref{sfq}). By~(\ref{granderoma}) we have $\dim|I_S(2)|=5$.
Let $\Sigma$ be the divisor on $|I_S(2)|$ of singular quadrics
(i.e.~of quadrics of rank at most $6$). Since $\dim|I_F(2)|=4$ and
every quadric containing $F$ is singular we have
\begin{equation}\label{split}
\Sigma=|I_F(2)|+\Sigma'
\end{equation}
where $\Sigma'$ is an effective divisor of degree $6$. The
hypersurface $Y$ is irreducible and non-degenerate and $Y\ss
supp(\Sigma)$, hence $Y\ss supp(\Sigma')$.
 By Item~(3) of Proposition~(\ref{explmap}) we know that
 $\deg(f_{\vv}\cl X\to Y)=2$ and
 by~(\ref{fujconst}) we have $\int_X f_{\vv}^{*}c_1(\cO_{\PP^5}(1))^4=12$, thus $\deg Y=6$:
 this gives that
 $\Sigma'=Y$ i.e.
 \begin{equation}
\Sigma=|I_F(2)|+Y.
 \end{equation}
\subsubsection{A closer view of $Y$}
Let $Fix(\phi_{\vv})\ss X$ be the locus of fixed points of
$\phi_{\vv}$; since $\phi_{\vv}$ is anti-symplectic
$Fix(\phi_{\vv})$ is a smooth Lagrangian surface. Let $\wh{X}\to
X$ be the blow up of $Fix(\phi_{\vv})$: then $\phi_{\vv}$ acts on
$\wh{X}$ with smooth quotient $\wh{Y}$. Since $\phi_{\vv}$ acts
trivially on $K_X$ the $4$-fold $\wh{Y}$ is a Calabi-Yau. The
natural map $\wh{Y}\to X/\la\phi_{\vv}\ra$ is a resolution of
singularities.
\begin{clm}\label{normale}
$Y$ is isomorphic to $X/\la\phi_{\vv}\ra$ and via this isomorphism
the map $f_{\vv}\cl X\to Y$ is identified with the quotient map
$X\to X/\la\phi_{\vv}\ra$.
\end{clm}
\begin{proof}
The map $f_{\vv}$ commutes with $\phi_{\vv}$ hence it descends to a map
$X/\la\phi_{\vv}\ra\to Y$ which is finite of degree $1$: we must show that this
map is an isomorphism. It suffices to show that $Y$ is normal, and since $Y$ is
a hypersurface this is equivalent to $Y$ being smooth in codimension $1$. The
Calabi-Yau $\wh{Y}$ is birational to $Y$ and hence any desingularization of $Y$
has Kodaira dimension equal to $0$; since $\deg Y=6$ we get by adjunction that
$Y$ is smooth in codimension $1$.
\end{proof}
It follows from the claim that $sing(Y)$ is a smooth surface and
that at a point $p\in sing(Y)$ the $4$-fold $Y$ is modelled on
$(\CC^2\tm V(x^2+y^2+z^2),0)$. The following result gives a
moduli-theoretic interpretation of the intersection $Y\cap
|I_F(2)|$.
\begin{clm}\label{goodrank}
We have $f_{\vv}^{*}|I_F(2)|=\Delta(\vv)+\Theta(\vv)$ where
$\Delta(\vv),\Theta(\vv)$ are given by~(\ref{deldef}) and
 Lemma~(\ref{delteta}) respectively.
\end{clm}
\begin{proof}
Let $\ww$ be as in~(\ref{wdef}); since $g=6$ the moduli space
$\cM(\ww)$ consists of a single point $[V]$ and $V$
satisfies~(\ref{sceltav}). Let $\b_V$ be the map
of~(\ref{betamap}). As is easily checked $\Im(\b_V)=|I_F(2)|$.
Proceeding as in the proof of Lemma~(\ref{nondegenere}), the case
$g=7$ (in the present proof $g=6$ but it makes no difference), we
get $\zeta_V\cl\PP(V^{\vee})\to\PP(H^0(V)^{\vee})$. For
$x\in\PP(V^{\vee})$ we let $F_x$ be the singular sheaf on $S$
defined in the proof of Lemma~(\ref{kappa1}); thus
$[F_x]\in\Delta(\vv)$. In fact we have an identification
\begin{equation}\label{isaraghi}
 \begin{matrix}
\PP(V^{\vee}) & \overset{\sim}{\lra} & \Delta(\vv)\\
x & \mapsto & [F_x]
 \end{matrix}
\end{equation}
Proposition~(\ref{explmap}) identifies $f_{\vv}$ with $q_{\vv}$
and hence~(\ref{cupiello}) becomes (via
Identification~(\ref{isaraghi})) the equality
$f_{\vv}|_{\Delta(\vv)}=\b_{V}\circ\zeta_V$; hence
$$f_{\vv}^{*}|I_F(2)|=\Delta(\vv)+\phi_{\vv}^{*}\Delta(\vv).$$
By~(\ref{phitet}) we have $\phi_{\vv}^{*}\Delta(\vv)=\Theta(\vv)$.
\end{proof}
\subsubsection{The dual of $Y$}
Going back to Strange duality we set $\ww:=1-\eta$ as in~(\ref{vedoppiav}). Thus
we have the map $f_{\ww}$ of~(\ref{fdoppiav}). Let $Y_{\ww}:=\Im(f_{\ww})$ and
$Y_{\vv}:=\Im(f_{\vv})$. Let $Y_{\vv}^{\vee}\ss |I_S(2)|^{\vee}$ be the dual of
$Y_{\vv}$; thus we have  a birational map

\begin{equation}\label{taumap}
\begin{matrix}
Y_{\vv}^{sm} & \overset{\tau_{\vv}}{\lra} &  Y_{\vv}^{\vee}\\
Q\ \ \ \  & \mapsto & \{Q'\in |I_S(2)|\ \ |\ sing(Q)\in Q'\}
\end{matrix}
\end{equation}
\begin{prp}\label{magicdual}
Keep notation as above. Then $Y_{\ww}=Y_{\vv}^{\vee}$.
\end{prp}
\begin{proof}
Let $Q\in f_{\vv}(U^{\flat}(\vv)\sm Fix(\phi_{\vv}))$. Then
$\rk(Q)=6$ and if $\{ p\}:=sing(Q)$ we have $p\notin S$. We claim
that there exists $[Z]\in S^{[2]}$ such that $p\in\la Z\ra$.
Assume the contrary. Let $\L\ss\PP^6$ be a a hyperplane not
containing $p$ and $Q_0:=Q\cap\L$. Projection from $p$ defines a
regular map $\pi\cl S\to Q_0$ which is is an embedding because  no
chord of $S$ contains $p$. The exact sequence of vector-bundles
$$0\to T_S\to\pi^{*}T_{Q_0}\to N_{S/Q_0}\to 0$$
gives that $\int_S c_2(N_{S/Q_0})=46$. Let $[\pi_{*}(S)]\in
H^4(Q_0)$ be the Poincar\`e dual of $\pi_{*}(S)$; since $\pi$ is
an embedding we have
$$\int_{Q_0}[\pi_{*}(S)]\wedge [\pi_{*}(S)]=\int_S
c_2(N_{S/Q_0})=46.$$
On the other hand we see directly that the left-hand side is equal
to $50$, contradiction. Thus there exists $[Z]\in S^{[2]}$ such
that $p\in\la Z\ra$. Clearly $\tau_{\vv}(p)=f_{\ww}([Z])$. Thus
$Y_{\vv}^{\vee}\ss Y_{\ww}$. Since both are irreducible
hypersurfaces we get that $Y_{\vv}^{\vee}\ss Y_{\ww}$.
\end{proof}
\begin{crl}
Keep notation as above. Then $\deg Y_{\ww}=6$ and $f_{\ww}$ has
degree $2$ onto its image.
\end{crl}
\begin{proof}
The map $f_{\ww}$ is base-point free, it commutes with the
involution $\phi_{\vv}$ and
$\int_{S^{[2]}}f_{\ww}^{*}c_1(\cO_{\PP^5}(1))^4=12$. Thus
$f_{\ww}$ is of finite degree $2d$ over its image $Y_{\ww}$ and
$\deg(Y_{\ww})=\frac{6}{d}$. On the other hand by
Proposition~(\ref{magicdual}) we know that any desingularization
of $Y_{\ww}$ has Kodaira dimension $0$ and hence by adjunction
$\deg(Y_{\ww})\ge 6$. Thus $\deg Y_{\ww}=6$ and $d=1$.
\end{proof}
\section{Connecting the examples}\label{tuttiperuno}
 \setcounter{equation}{0}
Let $\vv$ be given by~(\ref{vsymmetric}), and
assume that Hypothesis~(\ref{picint})
and~(\ref{eldpos}) hold. Assume also that
\begin{equation}\label{visa}
2n-2:=\la\vv,\vv\ra\le 4r-2.
\end{equation}
Let $h_{\vv}:=\theta_{\vv}(\eta-1)$ and let $H_{\vv}$ be a divisor on $\cM(\vv)$
such that $c_1(H_{\vv})=h_{\vv}$. Then $(\cM(\vv),H_{\vv})$ is a degree-$2$
polarized deformation of $(K3)^{[n]}$, see Item~(2) of Corollary~(\ref{amplelb})
for $r\ge 2$ and Subsubsection~(\ref{specifico}) for the case $r=1$. Choosing
different $\vv$'s we get many different families of degree-$2$ polarized
varieties of the same dimension: the methods that prove Theorem~(\ref{yoshi})
should also show that these varieties are {\sl polarized deformation
equivalent}, i.e.~that they are \lq\lq parametrized\rq\rq by the same connected
component of $\cQ_n^0$ (see~(\ref{spiego})). In other words we expect that given
$(\cM(\vv),H_{\vv})$ and $(\cM(\ww),H_{\ww})$ as above of the same dimension
there exist a proper submersive map of connected complex manifolds $\pi\cl\cX\to
B$, a relatively ample divisor $\cH$ on $\cX$, and $t,u\in B$ such that
$(X_t,H_t)\cong (\cM(\vv),H_{\vv})$ and $(X_u,H_u)\cong(\cM(\ww),H_{\ww})$. We
will prove that this is indeed the case in one significant example.  Let $(S,D)$
be a polarized $K3$ of degree $10$ and let
\begin{equation}\label{vudi}
\vv:=2+c_1(D)+2\eta.
\end{equation}
We assume that Hypothesis~(\ref{picint}) and~(\ref{eldpos}) hold with $\ell$
replaced by $c_1(D)$: thus $\cM(\vv)$ (stability is with respect to $D$)  is a
deformation of $(K3)^{[2]}$ and  $H_{\vv}$ is a degree-$2$ polarization of
$\cM(\vv)$. Let $T\ss\PP^3$ be a quartic surface not containing lines and let
$A$ be the (hyper)plane class on $T$. We let
$$\ww:=1+c_1(A)+\eta.$$
Then $\cM(\ww)=T^{[2]}$ and $H_{\ww}$ is ample of degree two.
\begin{prp}\label{sedeforma}
Keep notation and assumptions as above. Then $(\cM(\vv),H_{\vv})$ is polarized
deformation equivalent to $(\cM(\ww),H_{\ww})$.
\end{prp}
\begin{proof}
By surjectivity of the period map for $K3$'s there exists a quartic
$T_0\ss\PP^3$ containing a line $R$ and such that $H^{1,1}_{\ZZ}(T_0)=\ZZ
c_1(R)\op\ZZ c_1(A_0)$ where $A_0$ is the (hyper)plane class. Let
$$\ww_0:=1+c_1(A_0)+\eta.$$
The divisor $D_0:=2A_0-R$ is ample of degree $10$ thus $(T_0,D_0)$ is a
degree-$10$ polarized $K3$. Let
$$\vv_0:=2+c_1(D_0)+2\eta$$
and let $\cM(\vv_0)$ be the moduli space with (semi)stability with respect to
$D_0$. As is easily checked $D_0$ is $\vv_0$-generic and hence by
Theorem~(\ref{yoshi}) we know that $\cM(\vv_0)$ is smooth. We will define a
birational map
\begin{equation}\label{gammamap}
\g\cl \cM(\ww_0)\cdots>\cM(\vv_0).
\end{equation}
Let $[Z]\in T_0^{[2]}=\cM(\ww_0)$: an easy computation gives that
$$\dim\Ext^1(I_Z(A_0),\cO_{T_0}(A_0-R))=
\begin{cases}
1 & \mbox{if $Z\not\subset R$,}\\
2 & \mbox{if $Z\subset R$.}
\end{cases}$$
Furthermore any non-trivial extension
$$0\to \cO_{T_0}(A_0-R)\to E\to I_Z(A_0)\to 0$$
 is $D_0$-slope-stable and $v(E)=\vv_0$. Since for $[Z]\notin
 R^{(2)}$ we have a non-trivial extension $E_Z$ as above unique up to isomorphism,
we get a well-defined regular map
 \begin{equation}\label{pregamma}
 \begin{matrix}
(\cM(\ww_0)\sm R^{(2)}) & \lra & \cM(\vv_0)\\
 [Z] & \mapsto & E_Z.
\end{matrix}
\end{equation}
Let $[F]\in \cM(\vv_0)$; then $\chi(F(R-A_0))=1$ and since by stability we have
$h^2(F(R-A_0))=0$ we get that $h^0(F(R-A_0))\ge 1$. One checks easily
that~(\ref{pregamma}) is an isomorphism onto the open dense subset of
$\cM(\vv_0)$ parametrizing sheaves $F$ such that $h^0(F)=1$: we
define~(\ref{gammamap}) to be the birational map that corresponds
to~(\ref{pregamma}). One easily shows that $\g$ is the flop
(see~\cite{mukaisympl}) of $R^{(2)}$. Furthermore - and this is the main point -
we have $\g^{*}H_{\vv_0}=H_{\ww_0}$. If $\gamma$ were regular we would be done;
since $\gamma$ is not regular we proceed as follows. Let $\cX\to B_{\ww_0}$ be a
representative for the deformation space of $(\cM(\ww_0),H_{\ww_0})$,
i.e.~deformations of $\cM(\ww_0)$ that \lq\lq keep $H_{\ww_0}$ of type
$(1,1)$\rq\rq. Similarly let $\cX'\to B_{\vv_0}$ be a representative for the
deformation space of $(\cM(\vv_0),H_{\vv_0})$. Thus there is a divisor $\cH$ on
$\cX$ such that for every $s\in B_{\ww_0}$ the pull-back of $\cH$ to $X_s$, call
it $H_s$, is of degree $2$ and of course $H_0\sim H_{\ww_0}$. Similarly we have
$\cH'$ on $\cX'$. Now let $\L=R^{(2)}\ss T_0^{[2]}=\cM(\ww_0)$ and let $\L'\ss
\cM(\ww_0)$ be the corresponding copy of $R^{(2)}$ - the indeterminacy locus of
$\g^{-1}$. Let $B_{\ww_0}(\L)\ss B_{\ww_0}$ be the locus parametrizing
deformations of $(\cM(\ww_0),H_{\ww_0})$ for which $\L$ deforms too, and define
similarly $B_{\vv_0}(\L')$. By a Theorem of Voisin~\cite{voi} each of these loci
is smooth of codimension $1$.
\begin{clm}\label{lastclaim}
If we shrink enough $B_{\ww_0}$ and $B_{\vv_0}$ around $0$ the following holds.
Let $s\in (B_{\ww_0}\sm B_{\ww_0}(\L'))$. There exists $u\in (B_{\vv_0}\sm
B_{\vv_0}(\L')$ such that $(X_s,H_s)\cong (X'_u,H'_u)$.
\end{clm}
\begin{proof}
Let $\G\ss\cX$ be the locus swept out by the $\PP^2$'s which are deformations of
$\L$ and define similarly $\G'\ss\cX'$; thus we have $\PP^2$-bundles $\G\to
B_{\ww_0}(\L)$ and $\G'\to B_{\ww_0}(\L')$. Let $\pi\cl\cY\to\cX$ be the blow up
of $\G$. Let $E$ be the exceptional divisor of $\pi$; thus $\pi$ gives a
$\PP^2$-bundle $E\to\G$. Following Huybrechts~\cite{oldhuy} we see that $E$ has
another $\PP^2$-fibration structure $E\to\G'$ and that one can contract $\cY$
along this fibration and get a smooth $\cX''$. We still have a map $\cX''\to
B_{\vv_0}$ which is submersive, and the divisor $\cH''$ - the transform of
$\cH$. If $s\notin B_{\ww_0}(\L)$ then $(X''_s,H''_s)\cong (X_s,H_s)$. If $s\in
B_{\ww_0}(\L)$ then $X''_s$ is the flop of $X_s$ with center $\L_s$ (the
deformation of $\L$) and $H''_s$ is the divisor corresponding to $H_s$ via the
flop; in particular $X''_0\cong \cM(\vv_0)$. Considering the period map of
$\cX''$ we get that $\cX''\to B_{\ww_0}$ is the deformation space of
$(\cM(\vv_0),H_{\vv_0})$. The claim follows immediately.
\end{proof}
The proposition follows immediately from the above claim. In fact let $s\in
B_{\ww_0}$ parametrize $(\cM(\ww),H_{\ww})$ - such an $s$ exists as long as $T$
is sufficiently close to $T_0$. Since $T$ does not contain lines $s\notin
B_{\ww_0}(\L)$; by the claim there exists $\hat{u}\in B_{\vv_0}$ such that
$(\cM(\ww),H_{\ww})\cong (X'_{\hat{u}},H'_{\hat{u}})$. Since the moduli space of
polarized $K3$'s of degree $10$ (or any other degree) is irreducible
$(\cM(\vv),H_{\vv})$ is parametrized by a point  $\ov{u}\in B_{\vv_0}$ (again we
want $(S,D)$ sufficiently close to $(T_0,D_0)$) i.e.~$(\cM(\vv),H_{\vv})\cong
(X'_{\ov{u}},H'_{\ov{u}})$. Since the locus of $u\in B_{\vv_0}$ such that $H_u$
is ample is Zariski open we get the proposition.
\end{proof}
\vskip 1cm
 \scriptsize{
Kieran G. O'Grady\\
Universit\`a di Roma ``La Sapienza",\\
Dipartimento di Matematica ``Guido Castelnuovo",\\
Piazzale Aldo Moro n.~5, 00185 Rome, Italy,\\
e-mail: {\tt ogrady@mat.uniroma1.it}. }

\end{document}